\documentclass[10pt, letterpaper]{amsart}
\usepackage{amsmath, amsthm, amssymb, mathrsfs}
\usepackage{mathbbol}
\usepackage{txfonts}
\usepackage{graphicx}
\usepackage{psfrag}
\usepackage{graphicx}
\usepackage{verbatim}

\usepackage[usenames]{color}

\newtheorem{thm}{Theorem}[section]
\newcommand{\bt}{\begin{thm}}
\newcommand{\et}{\end{thm}}

\newtheorem{cor}[thm]{Corollary}
\newcommand{\bc}{\begin{cor}}
\newcommand{\ec}{\end{cor}}

\newtheorem{lem}[thm]{Lemma}
\newcommand{\bl}{\begin{lem}}
\newcommand{\el}{\end{lem}}

\newtheorem{prop}[thm]{Proposition}
\newcommand{\bp}{\begin{prop}}
\newcommand{\ep}{\end{prop}}

\newtheorem{defn}[thm]{Definition}
\newcommand{\bd}{\begin{defn}}      
\newcommand{\ed}{\end{defn}}

\newtheorem{rmrk}[thm]{Remark}
\newcommand{\br}{\begin{rmrk}}
\newcommand{\er}{\end{rmrk}}

\newtheorem{example}[thm]{Example}

\newcommand{\thmref}[1]{Theorem~\ref{#1}}
\newcommand{\secref}[1]{Section~\ref{#1}}
\newcommand{\lemref}[1]{Lemma~\ref{#1}}
\newcommand{\defref}[1]{Definition~\ref{#1}}
\newcommand{\corref}[1]{Corollary~\ref{#1}}
\newcommand{\propref}[1]{Proposition~\ref{#1}}
\newcommand{\remref}[1]{Remark~\ref{#1}}
\newcommand{\exref}[1]{Example~\ref{#1}}

\newcommand{\N}{\mathbb{N}}

\newcommand{\R}{\mathbb{R}}
\newcommand{\Z}{\mathbb{Z}}

\newcommand{\diam}{\operatorname{diam}}

\newcommand{\hm}{{\mathscr H}}
\newcommand{\lm}{{\mathscr L}}

\newcommand{\image}{\operatorname{Im}}

\newcommand{\Ricci}{\operatorname{Ric}}
\newcommand{\regset}{\mathcal{R}}

\newcommand{\lip}{\operatorname{Lip}}
\newcommand{\mass}[2][]{{\mathbf M_{#1}}(#2)}

\newcommand{\form}{{\mathcal D}}        
\newcommand{\curr}{{\mathbf M}}         
\newcommand{\intcurr}{{\mathbf I}}      

\newcommand{\Vol}{\operatorname{Vol}}

\newcommand{\fillvol}{{\operatorname{Fillvol}}}
\newcommand{\fillrad}{{\operatorname{Fillrad}}}

\newcommand{\rstr}{\:\mbox{\rule{0.1ex}{1.2ex}\rule{1.1ex}{0.1ex}}\:}
\newcommand{\bdry}{\partial}

\newcommand{\on}[1]{|_{#1}}
\newcommand{\spt}{\operatorname{spt}}
\newcommand{\ohne}{\backslash}

\newcommand{\injrad}{{\operatorname{inj}}}

\begin{document}

\title[Weak convergence and cancellation]{Weak convergence and cancellation}

\author[C.~Sormani]{Christina Sormani}
\thanks{C.~Sormani was partially supported by a PSC CUNY Research Grant. R.~Schul was partially supported by NSF grants DMS 0502747 and DMS 0800837. S.~Wenger was partially supported by NSF grant DMS 0707009.}
\address{C.~Sormani: CUNY Graduate Center and Lehman College}
\email{sormanic@member.ams.org}

\author[S.~Wenger]{Stefan Wenger}
\address
  {S.~Wenger: Department of Mathematics\\
University of Illinois at Chicago\\
851 S. Morgan Street\\
Chicago, IL 60607--7045 }
\email{wenger@math.uic.edu}

\author[]{an appendix by Raanan Schul and S.~Wenger}
\address{R.~Schul: UCLA Mathematics Department\\
Box 951555\\
Los Angeles, CA 90095-1555\\}
\email{schul@math.ucla.edu}

\date{February 16, 2009}


\begin{abstract}
 In this article, we study the relationship between the weak limit of a sequence of integral currents in a metric space and the possible Hausdorff limit of the sequence of supports. Due to cancellation, the weak limit is in general supported in a strict subset of the Hausdorff limit. We exhibit sufficient conditions in terms of topology of the supports which ensure that no cancellation occurs and that the support of the weak limit agrees with the Hausdorff limit of the supports. We use our results to prove countable $\hm^m$-rectifiability of the Gromov-Hausdorff limit of sequences of Lipschitz manifolds $M_n$ all of which are $\lambda$-linearly locally contractible up to some scale $r_0$. 
In the appendix, we show that the Gromov-Hausdorff limit need not be countably $\hm^m$-rectifiable if the $M_n$ have a common local geometric contractibility function which is only concave (and not linear). We also relate our results to work of Cheeger-Colding on the limits of noncollapsing sequences of manifolds with nonnegative Ricci curvature.
\end{abstract}

\maketitle

\section{Introduction and statement of the main results}

This article is concerned with the relationship between two notions of convergence, Hausdorff convergence for sets on the one hand, and weak convergence for (integral) currents on the other hand. These play an important role in problems in metric and Riemannian geometry and geometric measure theory, respectively. We give sufficient conditions in terms of topology on the sets (supports of the currents) that guarantee that these two notions `agree'. 

We work in the generality of metric spaces (for reasons that will become clear later).
In this setting, Ambrosio-Kirchheim have recently developed in \cite{Ambr-Kirch-curr} a powerful theory of metric integral currents which extends the classical Federer-Fleming theory \cite{Federer-Fleming} from Euclidean to (general) metric spaces $Z$. For $m\geq 0$, the space of metric integral $m$-currents in $Z$ is denoted by $\intcurr_m(Z)$. Similarly to the classical theory, elements $T\in\intcurr_m(Z)$ are `functionals' (on suitably generalized $m$-forms on $Z$) satisfying certain properties. 
As in the classical theory, there is a notion of mass, $\mass{T}$, an associated measure, $\|T\|$, and, for dimension $m\geq 1$, a
notion of boundary, $\bdry T$, which is an element of $\intcurr_{m-1}(Z)$. The support of $T$ is the closed set 
\begin{equation*}
 \spt T:= \{z\in Z: \text{ $\|T\|(B(z,\varepsilon))>0$ for all $\varepsilon>0$}\},
\end{equation*} 
 
where $B(z, \varepsilon)$ denotes the open ball of radius $\varepsilon$ centered at $z$.
Weak convergence is defined as pointwise convergence, as in the classical theory: $T_n$ converges weakly to $T$ if $T_n(fd\pi) \to T(fd\pi)$ for every Lipschitz $m$-form $fd\pi$ on $Z$. We refer to \secref{Section:Preliminaries} for definitions and details concerning integral currents.

We now turn to the relationship between weak convergence of the $T_n$ and Hausdorff convergence of the sets $\spt T_n$. For the definition of Hausdorff convergence of sets we refer to \secref{Section:Nagata}.
We say that  $(T_n)\subset\intcurr_m(Z)$ is a bounded sequence if
\begin{equation*}
 \sup_n[\mass{T_n} + \mass{\bdry T_n}]<\infty.
\end{equation*}
It is well-known, see e.g.~\lemref{lemma:incl-spt-limit}, that if $T_n$ converges weakly to some $T\in\intcurr_m(Z)$ and the sequence $(\spt T_n)$ converges in the Hausdorff sense to a closed subset $X\subset Z$ then $$\spt T\subseteq X.$$ 
The inclusion may be strict as the simple Examples \ref{example:torus} -- \ref{example:cylinder-gluing} illustrate. Strict inclusion is sometimes referred to as `cancellation', see \secref{Section:examples-cancellation}. In our first two theorems below we show that under suitable restrictions on the local topology of $\spt T_n$ cancellation cannot occur, thus $\spt T = X$. This is then used to prove metric structure results for the Gromov-Hausdorff limit of sequences of Riemannian manifolds which are uniformly linearly locally contractible up to some scale and sequences of noncollapsing Riemannian manifolds of nonnegative Ricci curvature. 

We turn to the precise statements of our results and first recall the following definition.

\bd
Given a metric space $Y$ and $\lambda\geq 1$, a subset $A\subset Y$ is called $\lambda$-linearly locally contractible  in $Y$ if, for all $a\in A$ and $r>0$, $A\cap \bar{B}(a, r)$ is contractible in $\bar{B}(a, \lambda r)$.  Given $m\geq 0$, a subset $A\subset Y$ is called $\lambda$-linearly locally $m$-connected in $Y$ if, for every $k\in\{0,\dots, m\}$, all $a\in A$ and $r>0$, every continuous map $f: S^k\to Y$ with image in $A\cap \bar{B}(a, r)$ has a continuous extension $\bar{f}: \bar{B}^{k+1}\to Y$ with image in $\bar{B}(a, \lambda r)$.
If for some $r_0>0$, every $r_0$-ball in $Y$ is $\lambda$-linearly locally contractible ($m$-connected) in $Y$ then $Y$ is called $\lambda$-linearly locally contractible ($m$-connected) up to scale $r_0$.

\ed
Here, $S^k$ and $\bar{B}^{k+1}$ denote the unit $k$-sphere and the closed unit ball in $\R^{k+1}$, respectively. We will also need the following definition, taken from \cite{Lang-Schlichenmaier}.

\bd\label{definition:lip-conn-small}
A metric space $Y$ is Lipschitz $m$-connected in the small if there exist constants $\gamma>0$ and $\delta>0$ such that for every $k\in\{0,\dots, m\}$, every $\nu$-Lipschitz map $f:S^k\to Y$ with $\nu\leq\delta$ has a $\gamma\nu$-Lipschitz extension $\bar{f}: \bar{B}^{k+1}\to Y$. 
\ed

Clearly, every compact Riemannian $n$-manifold is Lipschitz $m$-connected in the small for every $m$; the same is true for every compact metric space locally biLipschitz homeomorphic to an open set in $\R^n$. 

In its most general form, the main result of our paper can be stated as follows. (For the convenience of the geometrically minded reader, we will include versions for Riemanian manifolds of all our theorems in the corollaries below.)

\bt\label{theorem:main}
 Let $Z$ be a complete metric space, $m\geq 1$, and $(T_n)\subset\intcurr_m(Z)$ a bounded sequence of integral currents weakly converging to some $T\in\intcurr_m(Z)$. Suppose that the support $\spt T_n$ of each $T_n$ is compact, Lipschitz $m$-connected in the small, and has Hausdorff and Nagata dimension $< m+1$.
Let furthermore $z\in Z$ and $\lambda\geq 1$, $r>0$. If there exists a sequence $z_n\to z$ with $z_n\in\spt T_n$ and such that  $B(z_n, r)\cap \spt \bdry T_n=\emptyset$ and $B(z_n, r)\cap \spt T_n$ is $\lambda$-linearly locally $m$-connected in $\spt T_n$ for all $n$ large enough, then
\begin{equation}\label{equation:vol-growth-main-thm}
 \|T\|(B(z, s))\geq C\lambda^{-m(m+1)}s^m\quad\text{ for all $0\leq s\leq 2^{-(m+6)}\lambda^{-(m+1)}r$,}
\end{equation}
for some $C>0$ depending only on $m$.
In particular, $z\in\spt T$ and $\|T\|$ has strictly positive lower $m$-dimensional density at $z$.
\et

For the definition and properties of the Nagata dimension (also called Assouad-Nagata dimension in the literature) we refer to \secref{Section:Nagata} and to \cite{Lang-Schlichenmaier}. We mention here that the Nagata dimension of a metric space is always an integer (unlike the Hausdorff dimension which can be any nonnegative number) and is at  least its topological dimension. Furthermore, the Nagata dimension of any compact $m$-dimensional Riemannian manifold is $m$. This is more generally true of any compact metric space locally homeomorphic to an open subset of $\R^m$. Thus,  \thmref{theorem:main} implies:

\bc\label{Corollary:special-case-main-thm}
 Let $(M_n)$ be a sequence of compact oriented $m$-dimensional Riemannian manifolds with uniformly bounded volume and volume of boundary. Let $p_n\in M_n$ and suppose that for some $\lambda\geq 1$ and $r>0$ and each $n$, the ball $B_{M_n}(p_n, r)$ is $\lambda$-linearly locally contractible in $M_n$ and does not intersect $\bdry M_n$. Suppose further that $Z$ is a metric space and $\varphi_n:M_n\hookrightarrow Z$ are isometric embeddings. If $\varphi_n(p_n)$ converges to some $z$ and if the sequence of currents $T_n:=\varphi_{n\#}\Lbrack M_n\Rbrack$ weakly converges to some $T\in\intcurr_m(Z)$ then
 \begin{equation*}
 \|T\|(B(z, s))\geq C\lambda^{-m(m+1)}s^m\quad\text{ for all $0\leq s\leq 2^{-(m+6)}\lambda^{-(m+1)}r$,}
\end{equation*}
for some $C>0$ depending only on $m$.
In particular, $z\in\spt T$ and $\|T\|$ has strictly positive lower $m$-dimensional density at $z$.
\ec 

For the definition of $\Lbrack M_n\Rbrack$ see \eqref{equation:def-current-mfd}. The term {\it isometric} in the corollary means distance preserving (when $M_n$ is viewed as a metric space). Throughout this paper we assume manifolds to be connected.
In \exref{Example:Nagata-counter} we show that the hypothesis in \thmref{theorem:main} that $\spt T_n$ has Nagata and Hausdorff dimension $< m+1$ cannot be replaced by $\leq m+1$, even when $Z=\R^{m+1}$.  Note also that in \thmref{theorem:main} the constants appearing in \defref{definition:lip-conn-small} for $\spt T_n$ are allowed to `degenerate' as $n\to\infty$. 

\thmref{theorem:main} can be used to prove the following result:

\bt\label{Theorem:spt-rect}
 Let $Z$ be a complete metric space, $m\geq 1$, and $(T_n)\subset\intcurr_m(Z)$ a bounded sequence of integral currents weakly converging to some $T\in\intcurr_m(Z)$, with $\bdry T_n=0$ for all $n$. Suppose that the support $\spt T_n$ of each $T_n$ is compact, Lipschitz $m$-connected in the small, and has Hausdorff and Nagata dimension $< m+1$. If for some $\lambda\geq 1$ and $r>0$, $\spt T_n$ is $\lambda$-linearly locally $m$-connected up to scale $r$ for all $n$ large enough, then $\spt T_n$ converges in the Hausdorff sense to $\spt T$ and 
\begin{equation}\label{equation:growth-spt}
 \|T\|(B(z, s))\geq C\lambda^{-m(m+1)}s^m
\end{equation}
 for all $z\in\spt T$ and $0\leq s\leq 2^{-(m+6)}\lambda^{-(m+1)}r$, where $C>0$ depends only on $m$. In particular, $\spt T$ is countably $\hm^m$-rectifiable.
\et

We refer to \secref{Section:Hausdorff-measure} for the definition of countable $\hm^m$-rectifiability.
Note that the support of an arbitrary integral current $T$ in $Z$ need not be countably $\hm^m$-rectifiable; indeed if $Z$ is for example $Z$ is a separable Banach space then an integral current $T$ in $Z$ can have $\spt T = Z$.
Clearly, \eqref{equation:growth-spt} implies that $\|T\|$ has strictly positive lower $m$-dimensional density at every point $z\in\spt T$. 
The last claim in the theorem follows from this together with Theorem 4.6 in \cite{Ambr-Kirch-curr}. \thmref{Theorem:spt-rect} clearly implies the following (also compare with \corref{Corollary:Riemannian-lin-loc}):

\bc
 Let $(M_n)$ be a sequence of closed oriented $m$-dimensional Riemannian manifolds with uniformly bounded volume and such that for some $\lambda\geq 1$ and $r>0$, each $M_n$ is $\lambda$-linearly locally contractible up to scale $r$. Suppose further that $Z$ is a metric space and $\varphi_n:M_n\hookrightarrow Z$ are isometric embeddings. If the sequence of currents $T_n:=\varphi_{n\#}\Lbrack M_n\Rbrack$ weakly converges to some $T\in\intcurr_m(Z)$ then $\varphi_n(M_n)$ Hausdorff converges to $\spt T $. Furthermore $\|T\|$ has strictly positive lower $m$-dimensional density at every $z\in\spt T$ and, in particular, $\spt T$ is countably $\hm^m$-rectifiable.
\ec

A slightly more general form of \thmref{Theorem:spt-rect} will be discussed in \remref{Remark:abstract-sequence-limits}. 
In \secref{Section:Ricci} we will moreover prove an analog of \thmref{Theorem:spt-rect} for sequences of oriented closed Riemannian manifolds with nonnegative Ricci curvature, with a uniform upper bound on diameter and strictly positive lower bound on volume. As is well-known, such sequences admit an isometric embedding into some compact metric space and therein give rise to a sequence of integral currents. However, they do not satisfy the hypotheses of \thmref{Theorem:spt-rect} in general, see \secref{Section:Ricci} for details and references. Nevertheless, \thmref{Theorem:conv-Ricci} shows that the support of the weak limit is countably $\hm^m$-rectifiable and coincides with the (Gromov-)Hausdorff limit of the sequence. Countable $\hm^m$-rectifiability was proved before by Cheeger-Colding in \cite{ChCo-PartIII};  \thmref{Theorem:conv-Ricci} gives a new perspective of this fact.

\thmref{theorem:main} can be used to prove the following theorem whose statement does not involve any currents.

\bt\label{Theorem:GH-limit}
 Let $m\geq 1$ and let $(X_n)$ be a sequence of compact metric spaces such that each $X_n$ is locally biLipschitz homeomorphic to the open unit ball $B\subset \R^m$. Suppose further that each $X_n$ is orientable and that the $X_n$ have uniformly bounded diameter and Hausdorff $m$-measure. If there exist $\lambda\geq 1$ and $r>0$ such that $X_n$ is $\lambda$-linearly locally contractible up to scale $r$ for all $n$ large enough,  then there exists a subsequence $X_{n_j}$ which converges in the Gromov-Hausdorff sense to a compact and countably $\hm^m$-rectifiable space $X$ with $0<\hm^m(X)<\infty$.
\et

Here, {\it $X_n$ orientable} means that there exist finitely many biLipschitz maps $\varphi_i: B\hookrightarrow X_n$ such that the sets $\varphi_i(B)$ cover $X_n$ and such that $\det(\nabla(\varphi_i^{-1}\circ\varphi_j))>0$ almost everywhere on 
$\varphi_j^{-1}(\varphi_i(B))$.
%
%
As a simple consequence we obtain the following result.

\bc\label{Corollary:Riemannian-lin-loc}
 Let $(M_n)$ be a sequence of closed, orientable $m$-dimensional Riemannian manifolds of uniformly bounded volume. If there exist $\lambda\geq 1$ and $r>0$ such that $M_n$ is $\lambda$-linearly locally contractible up to scale $r$ for all $n$ large enough, then there exists a subsequence $(M_{n_j})$ which converges in the Gromov-Hausdorff sense to a compact and countably $\hm^m$-rectifiable metric space $X$ with $0<\hm^m(X)<\infty$.
\ec

 Note that no diameter bound is assumed in the corollary. The existence of a Gromov-Hausdorff limit in \corref{Corollary:Riemannian-lin-loc} was already proved by Greene-Petersen in \cite{Greene-Petersen} under the following weaker conditions: Let $r>0$ and let $\varrho : [0, r)\to [0,\infty)$ be a continuous function with $\varrho(0)=0$ and $\varrho(s)\geq s$ for all $s\in[0,r)$. Then $\varrho$ is called a local geometric contractibility function for a given metric space $X$ if for every $z\in X$ and every $
s\in(0,r)$ the ball $B(z,s)$ is contractible in $B(z, \varrho(s))$. It is shown in \cite{Greene-Petersen} that a sequence of closed $m$-dimensional Riemannian manifolds with a uniform upper bound on volume and such that $\varrho$ is a local geometric contractibility function for every $M_n$ then a subsequence $M_{n_j}$ converges in the Gromov-Hausdorff sense to a compact metric space $X$. Their proof relies on a lower bound -- proved in \cite{Greene-Petersen} using an argument similar to one in \cite{Gromov-filling} -- for the filling radius of (smooth approximations to) distance spheres in such Riemannian manifolds. In \thmref{Theorem:lower-bounds-fillrad} we prove a generalization of their filling radius estimate using a different approach. \thmref{Theorem:lower-bounds-fillrad} will be used in the proof of our main theorem.
Note that our assumption in \corref{Corollary:Riemannian-lin-loc} is that all $M_n$ have local geometric contractibility function $\varrho(s)= \lambda s$. 
In the appendix of the present paper, Raanan Schul and the second author show that if in \corref{Corollary:Riemannian-lin-loc} the condition on $\lambda$-linear local $m$-connectedness is replaced by uniform geometric contractibility,  
then the Gromov-Hausdorff limit $X$ need not be countably $\hm^m$-rectifiable. If, in addition, no uniform bound on volume is assumed, then $M_n$ can converge to an infinite dimensional space, as was shown by Ferry-Okun in \cite{Ferry-Okun}.

Our initial motivation for studying the relationship between Hausdorff convergence and weak convergence originates from the article \cite{Sormani-Wenger-long} in which we first introduced and studied a new distance, called flat intrinsic distance, between compact oriented Riemannian manifolds. Roughly, the flat intrinsic distance between $M$ and $M'$ is the infimal flat distance that can be achieved by isometrically embedding $M$ and $M'$ in a metric space $Z$. Flat intrinsic convergence of a sequence of Riemannian manifolds then amounts to flat convergence in a suitable metric space into which the sequence isometrically embeds. It is natural to ask what the relationship between the flat intrinsic and Gromov-Hausdorff convergence is. Since flat convergence implies weak convergence, our results above can be interpreted as giving sufficient conditions which ensure that a Gromov-Hausdorff limit coincides with the flat intrinsic limit, see \cite{Sormani-Wenger-long} and also \remref{Remark:abstract-sequence-limits} of the present paper. 

\bigskip

The paper is organized as follows. In \secref{Section:Lip-extension} we prove a Lipschitz extension theorem, \thmref{Theorem:Lip-extension-important}, which is a variant of a result of Lang-Schlichenmaier \cite{Lang-Schlichenmaier}. In \secref{Section:fillrad-bound} we use this theorem to exhibit lower bounds for the absolute filling radius of slices with spheres of integral currents, generalizing a theorem in  \cite{Greene-Petersen}.
\secref{Section:main-results} is devoted to the proofs of the results stated in the introduction. In \secref{Section:Ricci} we show that the (Gromov-) Hausdorff limit of sequences of Riemannian manifolds with nonnegative Ricci curvature agrees with the support of the weak limit and is countably $\hm^m$-rectifiable. The content of the appendix was discussed above.

\bigskip

{\bf Acknowledgments:} Some of the results of this paper were obtained while the second author was a Courant Instructor at New York University's Courant Institute of Mathematical Sciences and the first author was visiting Courant Institute. They wish to thank Courant Institute for the superb working environment.
\section{Preliminaries}\label{Section:Preliminaries}

\subsection{Hausdorff and Gromov-Hausdorff distance, Nagata dimension}\label{Section:Nagata}
Let $(Z,d)$ be a metric space and $A,B\subset Z$ closed. The Hausdorff distance between $A$ and $B$ in $Z$ is 
\begin{equation*}
 d_H(A,B):= \inf\left\{\varepsilon\geq 0: \text{$A\subset U_\varepsilon(B)$ and $B\subset U_\varepsilon(A)$}\right\},
\end{equation*}
where $U_\varepsilon(B)$ denotes the open $\varepsilon$-neighborhood of $B$. It is clear that $d_H$ is symmetric and satisfies the triangle inequality. If $d_H(A, B)=0$ then $A=B$.

\bd
A sequence $(A_n)$ of closed subsets of $Z$ is said to converge in the Hausdorff sense to a closed subset $A\subset Z$ if $d_H(A_n, A)\to0$ as $n\to\infty$.
\ed

If $X$ and $Y$ are metric spaces then the Gromov-Hausdorff distance between $X$ and $Y$ is
\begin{equation*}
 d_{GH}(X,Y):= \inf\left\{d_H(\iota_X(X), \iota_Y(Y)): \text{$Z$ metric space, $\iota_X: X\hookrightarrow Z$, $\iota_Y: Y\hookrightarrow Z$ isometric}\right\}.
\end{equation*}

A sequence of complete metric spaces $X_n$ is said to converge in the Gromov-Hausdorff sense to a complete metric space $X$ if $d_{GH}(X_n, X)\to 0$ as $n\to\infty$. By Gromov's famous compactness theorem for metric spaces, see \cite{Gromov-polynomial}, every uniformly compact sequence $(X_n)$ of metric spaces admits an isometric embedding into a common compact metric space $Z$. It follows that a subsequence $X_{n_j}$ converges in the Gromov-Hausdorff sense to a compact metric space, in fact to a closed subset of $Z$. Here, a sequence $(X_n)$ of compact metric spaces $X_n$ is called uniformly compact if $\sup\diam X_n<\infty$ and if for every $\varepsilon>0$ there exists $N(\varepsilon)$ such that every $X_n$ can be covered by at most $N(\varepsilon)$ balls of radius $\varepsilon$.

We turn to the definition and properties of the Assouad-Nagata dimension for a metric space $Y$. For a detailed account we refer to \cite{Lang-Schlichenmaier}.  A family $(B_i)_{i\in I}$ of subsets of $Y$ is called $D$-bounded if $\diam B_i\leq D$ for all $i\in I$. For $s\geq 0$, the $s$-multiplicity of the family is the infimum of all $k\geq 0$ such that every subset of $X$ with diameter $\leq s$ intersects at most $k$ members of the family. 

\bd
 The (Assouad-)Nagata dimension of a metric space $Y$ is the infimum of all integers $n$ with the following properties: there exists a constant $c>0$ such that for all $s>0$, $Y$ has a $cs$-bounded covering with $s$-multiplicity at most $n+1$.
\ed

It can be shown that the Nagata dimension of $Y$ is at least its topological dimension. If $Y=Y_1\cup Y_2$ then the Nagata dimension of $Y$ is the maximum of the Nagata dimensions of $Y_1$ and $Y_2$. Every subset of $\R^n$ with nonempty interior has Nagata dimension $n$. The Nagata dimension is invariant under biLipschitz homeomorphisms. It thus follows that every compact metric space locally biLipschitz homeomorphic to an open subset of $\R^n$ has Nagata dimension $n$. In particular, every compact Riemannian $n$-manifold has Nagata dimension $n$. For proofs of all these statements and many more properties see \cite{Lang-Schlichenmaier}.

\subsection{Hausdorff measure and countable $\hm^m$-rectifiability}\label{Section:Hausdorff-measure}
Let $Z$ be a metric space and $A\subset Z$. The Hausdorff $m$-dimensional measure of $A$ is 
\begin{equation*}
 \hm^m(A):= \lim_{\delta\searrow 0}\inf\left\{\sum_{i=1}^\infty \omega_m\left(\frac{\diam(B_i)}{2}\right)^m :
      B\subset \bigcup_{i=1}^\infty B_i\text{, }\diam(B_i)<\delta\right\},
\end{equation*}
where $\omega_m$ denotes the Lebesgue measure of the unit ball in $\R^m$. For $Z=\R^m$, $\hm^m$ agrees with the Lebesgue measure. The Hausdorff measure $\hm^s$ can also be defined for non-integer numbers $s$ by replacing the power $m$  in the above definition by $s$ and $\omega_m$ by a suitable positive number $\omega_s$, see for example \cite{Kirchheim}. The Haudorff dimension of $A$ is the infimum over all $s\geq 0$ such that $\hm^s(A)=0$.

\bd
An $\hm^m$-measurable set $A\subset Z$
is said to be countably $\hm^m$-rectifiable if there exist countably many Lipschitz maps $f_i :B_i\longrightarrow Z$ from subsets
$B_i\subset \R^m$ such that
\begin{equation*}
\hm^m\left(A\ohne \bigcup f_i(B_i)\right)=0.
\end{equation*}
\ed

Finally, the $m$-dimensional 
lower density $\Theta_{*m}(\mu, x)$ of a finite Borel measure $\mu$ on $Z$ at a point $z$ is given by the formula
\begin{equation*}
 \Theta_{*m}(\mu, z):= \liminf_{r\searrow 0}\frac{\mu(B(z,r))}{\omega_m r^m}.
\end{equation*}

\subsection{Currents in metric spaces}\label{section:currents}
We recall the basic definitions from the theory of currents which we need in this paper. The general reference is \cite{Ambr-Kirch-curr}.

Let $(Z,d)$ be a complete metric space and $m\geq 0$ and let $\form^m(Z)$ be the set of $(m+1)$-tuples $(f,\pi_1,\dots,\pi_m)$ 
of Lipschitz functions on $Z$ with $f$ bounded. Sometimes, we write the short-hand $fd\pi$ for $(f, \pi_1,\dots, \pi_m)$. The Lipschitz constant of a Lipschitz function $f$ on $Z$ will
be denoted by $\lip(f)$.
\bd
An $m$-dimensional metric current  $T$ on $Z$ is a multi-linear functional on $\form^m(Z)$ satisfying the following
properties:
\begin{enumerate}
 \item If $\pi^j_i$ converges pointwise to $\pi_i$ as $j\to\infty$ and if $\sup_{i,j}\lip(\pi^j_i)<\infty$ then
       \begin{equation*}
         T(f,\pi^j_1,\dots,\pi^j_m) \longrightarrow T(f,\pi_1,\dots,\pi_m).
       \end{equation*}
 \item If $\{z\in Z:f(z)\not=0\}$ is contained in the union $\bigcup_{i=1}^mB_i$ of Borel sets $B_i$ and if $\pi_i$ is constant 
       on $B_i$ then
       \begin{equation*}
         T(f,\pi_1,\dots,\pi_m)=0.
       \end{equation*}
 \item There exists a finite Borel measure $\mu$ on $Z$ such that
       \begin{equation}\label{equation:mass-def}
        |T(f,\pi_1,\dots,\pi_m)|\leq \prod_{i=1}^m\lip(\pi_i)\int_Z|f|d\mu
       \end{equation}
       for all $(f,\pi_1,\dots,\pi_m)\in\form^m(Z)$.
\end{enumerate}
\ed
The space of $m$-dimensional metric currents on $Z$ is denoted by $\curr_m(Z)$ and the minimal Borel measure $\mu$
satisfying \eqref{equation:mass-def} is called mass of $T$ and written as $\|T\|$. We also call mass of $T$ the number $\|T\|(Z)$ 
which we denote by $\mass{T}$.
The support of $T$ is, by definition, the closed set $\spt T$ of points $z\in Z$ such that $\|T\|(B(z,r))>0$ for all $r>0$. 

Every function $\theta\in L^1(K,\R)$ with $K\subset\R^m$ Borel measurable induces an element of $\curr_m(\R^m)$ by
\begin{equation*}
 \Lbrack\theta\Rbrack(f,\pi_1,\dots,\pi_m):=\int_K\theta f\det\left(\frac{\partial\pi_i}{\partial x_j}\right)\,d\lm^m
\end{equation*}
for all $(f,\pi_1,\dots,\pi_m)\in\form^m(\R^m)$.

The restriction of $T\in\curr_m(Z)$ to a Borel set $A\subset Z$ is given by 
\begin{equation*}
  (T\rstr A)(f,\pi_1,\dots,\pi_m):= T(f\chi_A,\pi_1,\dots,\pi_m).
\end{equation*}
This expression is well-defined since $T$ can be extended to a functional on tuples for which the first argument lies in 
$L^\infty(Z,\|T\|)$. It can be shown that $T(f,\pi_1, \dots, \pi_m)$ only depends on the values of $f, \pi_1, \dots, \pi_m$ on $\spt T$.

If $m\geq 1$ and $T\in\curr_m(Z)$ then the boundary of $T$ is the functional
\begin{equation*}
 \bdry T(f,\pi_1,\dots,\pi_{m-1}):= T(1,f,\pi_1,\dots,\pi_{m-1}).
\end{equation*}
It is clear that $\bdry T$ satisfies conditions (i) and (ii) in the above definition. If $\bdry T$ also satisfies (iii) then $T$ is called a normal current.
By convention, elements of $\curr_0(Z)$ are also called normal currents.

The push-forward of $T\in\curr_m(Z)$ 
under a Lipschitz map $\varphi$ from $Z$ to another complete metric space $Y$ is given by
\begin{equation*}
 \varphi_\# T(g,\tau_1,\dots,\tau_m):= T(g\circ\varphi, \tau_1\circ\varphi,\dots,\tau_m\circ\varphi)
\end{equation*}
for $(g,\tau_1,\dots,\tau_k)\in\form^m(Y)$. This defines an $m$-dimensional metric current on $Y$.
It follows directly from the definitions that $\bdry(\varphi_{\#}T) = \varphi_{\#}(\bdry T)$.

\bd
 A sequence $(T_n)$ in $\curr_m(Z)$ is said to converge weakly to $T\in\curr_m(Z)$ if $$T_n(f, \pi_1,\dots,\pi_m)\to T(f,\pi_1,\dots,\pi_m)$$ for all $(f,\pi_1,\dots,\pi_m)\in\form^m(Z)$.
\ed

In this article we will exclusively work with integral currents, defined below: 
An element $T\in\curr_0(Z)$ is called integer rectifiable if there exist finitely many points $x_1,\dots,x_n\in X$ and $\theta_1,\dots,\theta_n\in\Z\ohne\{0\}$ such
that
\begin{equation*}
 T(f)=\sum_{i=1}^n\theta_if(x_i)
\end{equation*}
for all bounded Lipschitz functions $f$.

\sloppy
 A current $T\in\curr_m(Z)$ with $m\geq 1$ is said to be integer rectifiable if the following properties hold:
 \begin{enumerate}
  \item $\|T\|$ is concentrated on a countably $\hm^m$-rectifiable set and vanishes on $\hm^m$-negligible Borel sets.
  \item For any Lipschitz map $\varphi:Z\to\R^m$ and any open set $U\subset Z$ there exists $\theta\in L^1(\R^m,\Z)$ such that 
    $\varphi_\#(T\rstr U)=\Lbrack\theta\Rbrack$.
 \end{enumerate}
 \fussy
Integer rectifiable normal currents are called integral currents. The corresponding space is denoted by $\intcurr_m(Z)$. If $A\subset\R^m$ is a Borel set of finite measure and
finite perimeter then $\Lbrack\chi_A\Rbrack \in\intcurr_m(\R^m)$. Here, $\chi_A$ denotes the characteristic function. If $T\in\intcurr_m(Z)$ and if $\varphi:Z\to Y$ is a Lipschitz 
map into another complete metric space then $\varphi_{\#}T\in\intcurr_k(Y)$. If $Z=\R^n$ then there is a natural isometric isomorphism between the space of compactly supported Ambrosio-Kirchheim integral (integer rectifiable) $m$-currents and the space of Federer-Fleming integral (integer rectifiable) $m$-currents in $\R^n$, see \cite{Ambr-Kirch-curr}.

Given a Riemannian manifold $N$ and an oriented $m$-dimensional submanifold $M$ of $N$, possibly with boundary. Let $\tau = \tau_1\wedge\dots\wedge\tau_m$ be the orientation of $M$, where $\{\tau_1, \dots,\tau_m\}$ is an $m$-tuple of orthonormal vector fields. If $M$ has finite volume and finite volume of the boundary then it gives rise to an integral current $\Lbrack M\Rbrack\in\intcurr_m(N)$ by
\begin{equation}\label{equation:def-current-mfd}
 \Lbrack M\Rbrack(f,\pi_1,\dots, \pi_m):= \int_Mf\,\langle d\pi_1\wedge\dots\wedge d\pi_m, \tau\rangle d\hm^m
\end{equation}
for all $(f,\pi_1,\dots, \pi_m)\in\form^m(N)$. By definition,
\begin{equation*}
 \langle d\pi_1\wedge\dots\wedge d\pi_m, \tau\rangle(z):= \det\left(d_z\pi_i(\tau_j(z))\right)
\end{equation*}
for all $z\in M$ at which each $\pi_i$ is differentiable (which is for $\hm^m$-almost every $z$ by Rademacher's theorem). Note that $\hm^m\rstr M$ agrees with the Riemannian volume on $M$.

\bigskip

We conclude this section with the following easy result.

\bl\label{lemma:incl-spt-limit}
 Let $(T_m)\subset\intcurr_m(Z)$ be a bounded sequence converging weakly to some $T\in\intcurr_m(Z)$. Then $\spt T$ is contained in the set of points $z\in Z$ for which there exists a sequence $z_n$ converging to $z$ with $z_n\in\spt T_n$ for all $n$ large enough. In particular, if $(\spt T_n)$ converges in the Hausdorff sense to a closed subset $X\subset Z$ then $\spt T\subset X$.
\el

\begin{proof}
 Let $z\in \spt T$ and $\varepsilon>0$. Since $\|T\|(B(z,\varepsilon/2))>0$, there exist $f, \pi_1,\dots, \pi_m\in\lip(Z)$ such that $\spt f\subset B(z, \varepsilon)$ and $T(f,\pi_1,\dots, \pi_m)\not=0$. It follows that $T_n(f,\pi_1,\dots, \pi_m)\not=0$ for $n$ large enough and we thus conclude from Proposition 2.7 in \cite{Ambr-Kirch-curr} that
$\|T_n\|(B(z, \varepsilon))>0$ and $\spt T_n\cap B(z, \varepsilon)\not=\emptyset$.
\end{proof}

\section{A Lipschitz extension theorem}\label{Section:Lip-extension}

The main purpose of this section is to prove \thmref{Theorem:Lip-extension-important}, a Lipschitz extension theorem which will be needed in the proof of \thmref{theorem:main}.  The results in this section are variants of results of Lang-Schlichenmaier in \cite{Lang-Schlichenmaier} and our proofs here follow closely those in \cite{Lang-Schlichenmaier}. 

\bd
Let $Y$ be a metric space, $n\geq 0$, and $\varrho\geq 1$. A subset $A\subset Y$ will be called $\varrho$-linearly locally weakly Lipschitz $n$-connected in $Y$ if for every $\nu'>0$ there exists $\gamma>0$ such that for all $k\in\{0,1,\dots, n\}$, all $a\in A$ and $r>0$, every $\nu$-Lipschitz map $f:S^k\to Y$ with  $\nu\geq\nu'$ and image in $A\cap \bar{B}(a, r)$ has a $\gamma\nu$-Lipschitz extension $\bar{f}: \bar{B}^{k+1}\to Y$ with image in $\bar{B}(a, \varrho r)$.
\ed

The following proposition is a variation of Theorem 5.1 in \cite{Lang-Schlichenmaier}.

\bp\label{proposition:weakly-Lip-connected}
Let $Y$ be a compact metric space, $\varrho\geq 1$, $r_0>0$, and $n\geq 0$. Suppose that $Y$ is Lipschitz $n$-connected in the small. If a ball $\bar{B}(y_0, r_0)$ in $Y$ is $\varrho$-linearly locally $n$-connected in $Y$ then $\bar{B}(y_0, r_0)$ is $(2\varrho)$-linearly locally weakly Lipschitz $n$-connected.
\ep

Up to a minor adjustment, the proof is the same as that in \cite{Lang-Schlichenmaier}. It consists of three parts, only the first of which needs modifying.

\begin{proof}
 Since $Y$ is Lipschitz $n$-connected in the small, there exist $\delta', \gamma'>0$ such that for all $l\in\{0,\dots, n\}$ every $\delta$-Lipschitz map $\varphi: \bdry [0,1]^{l+1}\to Y$ with $\delta\leq\delta'$ has a $\gamma'\delta$-Lipschitz extension to $[0,1]^{l+1}$.
 Fix $l\in\{0,\dots, n\}$ and put $D:= [0,1]^{l+1}$.
Let $f:\bdry D\to Y$ be a $\lambda$-Lipschitz map with image in $\bar{B}(y, r)\cap \bar{B}(y_0, r_0)$ for some $y\in\bar{B}(y_0, r_0)$ and some $r\leq 2r_0$. We show that $f$ has a Lipschitz extension $\bar{f}: D\to Y$ with image in $\bar{B}(y, \frac{3}{2}\varrho r)$. For this, pick a continuous extension $g: D\to Y$ of $f$ with image in $\bar{B}(y, \varrho r)$. This is possible because $\bar{B}(y_0, r_0)$ is $\varrho$-linearly locally $n$-connected. Let $\varepsilon>0$ be small enough, to be determined later. Equip $D$ with the structure of a piecewise Euclidean polyhedral complex whose $0$-skeleton $D^{(0)}$ is $D\cap ((1/N)\Z)^{l+1}$ and whose $(l+1)$-cells are cubes of edge length $1/N$, where $N\geq 2$ is chosen large enough so that $\lambda\leq \varepsilon N$ and $\diam g(C)\leq \varepsilon$ for every $(l+1)$-cell $C$ of $D$. For $k=0, \dots, n$, we successively find extensions $f^{(k)}: \bdry D \cup D^{(k)}\to Y$ of $f$ such that
$f^{(k)}\on{\bdry C}$ is $(2\gamma')^k\varepsilon N$-Lipschitz on every $(k+1)$-cell $C$ not contained in $\bdry D$ and such that 
\begin{equation*}
 f^{(k)}\left(\bdry D\cup D^{(k)}\right) \subset \bar{B}\left(y, \varrho r + 2^{-1}\varepsilon\gamma'\tau(k)\right)
\end{equation*}
where $\tau(k) = \frac{(2\gamma')^k-1}{2\gamma' - 1}$. For $k=0$, the extension $f^{(0)}:= g\on{\bdry D\cup D^{(0)}}$ satisfies these properties. Now, suppose that for a given $k$, $f^{(k)}$ is an extension with these properties. If $(2\gamma')^k\varepsilon\leq \delta'$ then there exists an extension $f^{(k+1)}: \bdry D\cup D^{(k+1)}\to Y$ of $f$ such that $f^{(k+1)}\on{C}$ is $\gamma'(2\gamma')^k\varepsilon N$-Lipschitz for every $(k+1)$-cell $C$. In particular, $f^{(k+1)}(C)$ is contained in the $2^{-1}\gamma'(2\gamma')^k\varepsilon$-neighborhood of $f^{(k)}(\bdry C)$. If $k+1\leq l$, it follows that for every $(k+2)$-cell $C'$ that is not contained in $\bdry D$, $f^{(k+1)}\on{\bdry C'}$ is $(2\gamma')^{k+1}\varepsilon N$-Lipschitz. 
If $\varepsilon>0$ is chosen sufficiently small so that, in particular, $\varepsilon\gamma'\tau(l+1)\leq \varrho r$, this iterative procedure gives the desired extension $\bar{f}:=f^{(l+1)}$ of $f$. This proves the first part. The second and third parts of the proof are the same as in \cite{Lang-Schlichenmaier}.
\end{proof}

The next theorem is an adjustment of Theorem 5.2 in \cite{Lang-Schlichenmaier} to our setting. We will indicate below the changes one has to make in Lang-Schlichenmaier's proof to obtain our theorem.

\bt\label{Theorem:LS-Thm5.2-gen}
 Let $X$ and $Y$ be metric spaces, and let $\alpha, \beta>0$ and $n\geq 1$. Suppose that $Z\subset X$ is a nonempty closed set and $(B_i)_{i\in I}$ is a covering of $X\backslash Z$ by subsets of $X\backslash Z$ such that
 \begin{enumerate}
  \item $\diam B_i \leq \alpha d(B_i, Z)$ for every $i\in I$,
  \item every set $D\subset X\backslash Z$ with $\diam D\leq \beta d(D, Z)$ meets at most $n+1$ members of $(B_i)_{i\in I}$. 
 \end{enumerate}
 Suppose that $Y$ is Lipschitz $(n-1)$-connected in the small and
 that $\bar{B}(y_0, r_0)$ is a ball in $Y$ which is $\varrho$-linearly locally weakly Lipschitz $(n-1)$-connected in $Y$. Then for all $\lambda\geq 1$ and $0<r<\frac{r_0}{16\lambda\varrho^n}$, every $\lambda$-Lipschitz map $f: Z\to Y$ with image in $\bar{B}(y_0, r_0/2)$ has a Lipschitz extension $\bar{f}: U_r(Z)\to Y$ with image in $U_{r'}(f(Z))$ where $r' = 8\lambda\varrho^n r$.
\et

Here, $U_r(Z)$ denotes the open $r$-neighborhood of $Z$ in $X$.
If $Z$ has Nagata dimension $\leq n-1$ then a covering $(B_i)_{i\in I}$ of $X\backslash Z$ satisfying the properties in the above theorem is constructed  in the proof of Theorem 1.6 in \cite{Lang-Schlichenmaier}. Therefore, we obtain the following result.

\bt\label{Theorem:Lip-extension-important}
 Let $X$ and $Y$ be metric spaces, let $n\geq 1$, $\varrho\geq 1$,  and $r_0>0$. Suppose that $Z\subset X$ is a nonempty closed set of Nagata dimension $\leq n-1$. Suppose further that $Y$ is Lipschitz $(n-1)$-connected in the small and
 that $\bar{B}(y_0,r_0)$ is a ball in $Y$ which is $\varrho$-linearly locally weakly Lipschitz $(n-1)$-connected in $Y$. Then for every $\lambda\geq 1$ and $0<r<\frac{r_0}{16\lambda\varrho^n}$, every $\lambda$-Lipschitz map $f: Z\to Y$ with image in $\bar{B}(y_0, r_0/2)$ has  a Lipschitz extension $\bar{f}: U_r(Z)\to Y$ with image in $U_{r'}(f(Z))$ where $r' = 8\lambda\varrho^n r$.
\et

\begin{proof}[Proof of \thmref{Theorem:LS-Thm5.2-gen}]
 The proof follows that of Theorem 5.2 in \cite{Lang-Schlichenmaier}. We need to make several modifications, however. Let $f: Z\to Y$ be a $\lambda$-Lipschitz map with image in $\bar{B}(y, r_0/2)$. Set $\delta:= \beta/(2(\beta+1))$ and, for each $i\in I$, let $\sigma_i: X\backslash Z\to \R$ be the $1$-Lipschitz function defined by
 \begin{equation*}
  \sigma_i(x):= \max\{0, \delta d(B_i, Z) - d(x, B_i)\}.
 \end{equation*} 
 It can then be shown (see \cite{Lang-Schlichenmaier}) that for every $x\in X\backslash Z$, there are at most $n+1$ indices $i\in I$ such that $\sigma_i(x)>0$. Now define $\bar{\sigma}:= \sum_{i\in I}\sigma_i$ and note that $\bar{\sigma}>0$ on $X\backslash Z$ because $(B_i)_{i\in I}$ covers $X\backslash Z$.
 Now define $g:X\backslash Z\to l^2(I)$ by
 \begin{equation*}
  g(x):= \left(\frac{\sigma_i(x)}{\bar{\sigma}(x)}\right)_{i\in I}
 \end{equation*}
 and note that the image of $g$ lies in the $n$-skeleton $\Sigma^{(n)}$ of the simplex $\Sigma:= \{(v_i)_{i\in I}: v_i\geq 0, \sum_{i\in I} v_i = 1\}\subset l^2(I)$. For every $i\in I$, choose a point $z_i\in Z$ such that
 \begin{equation*}
  d(z_i, B_i)\leq (2-\delta)d(B_i, Z).
 \end{equation*}
 Let $h^{(0)}: \Sigma^{(0)}\to Y$ be the map defined by $h^{(0)}(e_i):= f(z_i)$, where $e_i$ is the $i$-th vertex of $\Sigma$. So far, the proof was exactly the same as the one in \cite{Lang-Schlichenmaier}. 
 
 Now, let $\hat{\Sigma}$ be the subcomplex of $\Sigma$ consisting of those simplices $[e_{i_0}, \dots, e_{i_k}]\subset\Sigma^{(k)}$ for which there exists $x\in U_r(Z)$ with $\sigma_{i_j}(x)>0$ for every $j=0, \dots, k$. Note that actually $g(U_r(Z)\backslash Z)$ lies in $\hat{\Sigma}$.
Denote by $h^{(0)}$ again the restriction to $\hat{\Sigma}^{(0)}$ of the map $h^{(0)}$ defined above. Since $Y$ is Lipschitz $(n-1)$-connected in the small there exist constants $\delta'\in (0, 1/2)$ and $\gamma'\geq 1$ such that for all $k\in\{0, \dots, n-1\}$, every $\nu$-Lipschitz map $\varphi: \bdry\Delta^{k+1}\to Y$ with $\nu\leq\delta'$ has a $\gamma'\nu$-Lipschitz extension $\bar{\varphi}: \Delta^{k+1}\to Y$. Set $\nu':= \delta'\gamma'^{-1}\min\{r,1\}$ and choose $\gamma\geq 1$ such that for all $k\in\{0,\dots, n-1\}$, every $\nu$-Lipschitz map $\psi: \bdry\Delta^{k+1}\to Y$ with $\nu\geq \nu'$ and image in $\bar{B}(y, s)\cap \bar{B}(y_0, r_0)$ for some $y\in\bar{B}(y_0, r_0)$ and $s\leq 2r_0$ has a $\gamma\nu$-Lipschitz extension $\bar{\psi}: \Delta^{k+1}\to Y$ with image in $\bar{B}(y, \varrho s)$.  This is possible because $\bar{B}(y_0, r_0)$ is $\varrho$-linearly locally weakly Lipschitz $(n-1)$-connected.
Next, let $[e_i, e_j]\subset\hat{\Sigma}^{(1)}$ be a $1$-simplex. There then exists $x\in U_r(Z)$ with $\sigma_i(x)>0$ and $\sigma_j(x)>0$ and thus
 \begin{equation*}
  \begin{split}
   d(z_i, z_j) &\leq d(z_i, B_i) + d(B_i, x) + d(x, B_j) + d(B_j, z_j)\\
    &\leq d(z_i, B_i) + \delta d(B_i, Z) + \delta d(B_j, Z) + d(B_j, z_j)\\
    & \leq 2[d(B_i, Z) + d(B_j, Z)]\\
    &\leq 8r.
  \end{split}
 \end{equation*}
 The last inequality follows from the fact that
 \begin{equation*}
  d(Z, B_i) \leq d(Z, x) + d(x, B_i) < r +\delta d(B_i, Z)
 \end{equation*}
 and that $(1-\delta)^{-1}<2$. 
 It follows that $d(h^{(0)}(e_i), h^{(0)}(e_j))\leq \lambda d(z_i, z_j)\leq 8r\lambda$. 
For $m=0, 1,\dots, n-1$, successively extend $h^{(m)}$ to $h^{(m+1)}: \hat{\Sigma}^{m+1}\to Y$ as follows. Let $S\subset\hat{\Sigma}^{(m+1)}$ be an $(m+1)$-simplex and suppose that $h^{(m)}(\bdry S)$ is contained in the $8\lambda\varrho^mr$-ball around $h^{(m)}(e_i)=f(z_i)$ for some $e_i\in S^{(0)}$. (Note that this is the case for $m=0$.) We distinguish two cases: If $\lip(h^{(m)}\on{\bdry S})< \nu'$ then
$h^{(m)}\on{\bdry S}$ has a $\gamma'\lip(h^{(m)}\on{\bdry S})$-Lipschitz extension $h^{(m+1)}\on{S}$ to $S$. In particular, $h^{(m+1)}(S)$ is contained in the $r$-ball around $h^{(m)}(e_i)$. If, on the other hand, $\lip(h^{(m)}\on{\bdry S})\geq \nu'$ then $h^{(m)}\on{\bdry S}$ has a $\gamma\lip(h^{(m)}\on{\bdry S})$-Lipschitz extension $h^{(m+1)}\on{S}$ to $S$ with image in the $8\lambda\varrho^{m+1}r$-ball around $h^{(m)}(e_i)=f(z_i)$. We set $h:= h^{(n)}: \hat{\Sigma}^{(n)}\to Y$ and note that $h$ is Lipschitz on every simplex $S$ of $\hat{\Sigma}^{(n)}$ with Lipschitz constant
\begin{equation*}
 \lip(h\on{S})\leq C_1 \lip(h\on{S^{(0)}}) = \frac{C_1 \diam h(S^{(0)})}{\sqrt{2}}
\end{equation*}
for some constant $C_1$ depending on $r$, $\lambda$, $\gamma$, $\gamma'$, $\delta'$, $\varrho$, and $n$. Furthermore, $h(\hat{\Sigma})$ is contained in $U_{r'}(f(Z))$ with $r'=8\lambda\varrho^{n}r$. Finally, we define the extension $\bar{f}: U_r(Z)\to Y$ of $f$ such that $\bar{f} = h\circ g$ on $U_r(Z)$. The Lipschitz property of $\bar{f}$ follows now exactly as in \cite{Lang-Schlichenmaier} after (5.6). This completes the proof of our theorem.
\end{proof}

\section{Filling radius estimates and linear local $m$-connectedness}\label{Section:fillrad-bound}

Before stating the main result of this section we recall the notion of filling radius (first introduced and studied by Gromov in \cite{Gromov-filling}). Let $Z$ be a complete metric space and let $l^\infty(Z)$ be the Banach space of bounded functions on $Z$ with the supremum norm. Fix $z_0\in Z$. Then the map $\iota: Z\hookrightarrow l^\infty(Z)$ given by $\iota(z):= d(z, \cdot) - d(z_0,\cdot)$ defines an isometric embedding, called Kuratowski embedding. Let now $m\geq 0$ and let $T\in\intcurr_m(Z)$ be an integral current. The absolute filling radius of $T$ is 
\begin{equation*}
 \fillrad_\infty(T) = \inf\left\{\varepsilon\geq0: \text{$\exists S\in\intcurr_{m+1}(l^\infty(Z))$ with $\bdry S= \iota_\#T$ and $\spt S\subset U_\varepsilon (\iota(\spt T))$}\right\},
\end{equation*} 
where $U_\varepsilon(\iota(\spt T))$ denotes the open $\varepsilon$-neighborhood of $\iota(\spt T)$ in $l^\infty(Z)$. It is not difficult to see that $\fillrad_\infty(T)$ remains unchanged if $l^\infty(Z)$ is replaced by any other injective metric space $Z'$ into which $\spt T$ isometrically embeds. This justifies the word ``absolute'' filling radius. The main theorem of the present section is:

\bt\label{Theorem:lower-bounds-fillrad}
 Let $Z$ be a complete metric space and $m\geq 1$. Let $T\in\intcurr_m(Z)$ and suppose $\spt T$ is compact, Lipschitz $m$-connected in the small, and has Hausdorff and Nagata dimension $< m+1$. Let $z\in\spt T$ and suppose that for some $\lambda\geq 1$ and $r_0>0$ the ball $\bar{B}(z, r_0)\cap\spt T$ is $\lambda$-linearly locally $m$-connected in $\spt T$. Then 
 \begin{equation*}
 \fillrad_\infty(\bdry(T\rstr \bar{B}(z,r)))\geq \frac{r}{8(2\lambda)^{m+1}}
\end{equation*}
for almost every $r\in \left[0, \min\{2^{-(m+6)}\lambda^{-(m+1)}r_0, d(z, \spt\bdry T)\}\right]$.
\et

In particular, it follows that
\begin{equation*}
 \|T\|(\bar{B}(x,r))\geq\fillvol_\infty(\bdry(T\rstr \bar{B}(z,r)))\geq \frac{C'r^m}{[8(2\lambda)^{m+1}]^m}
\end{equation*}
for some $C'>0$ depending only on $m$, see \cite{Gromov-filling} or \cite{Ambr-Kirch-curr}, \cite{Wenger-GAFA}.

In the setting of Riemannian manifolds, \thmref{Theorem:lower-bounds-fillrad} was essentially proved by Greene-Petersen in \cite{Greene-Petersen}. Indeed, let $M$ be an $m$-dimensional Riemannian manifold, $z\in M$, and let $f_\varepsilon$ be smooth approximations of the distance function to $z$ on $M$. Then for almost every $r$ the subset $B_\varepsilon(r):=\{f_\varepsilon \leq r\}$ has smooth boundary. Greene-Petersen use an extension argument similar to the one in the proof of Gromov's Lemma 1.2.B in \cite{Gromov-filling} to show that if $M$ has a local geometric contractibility function $\varrho: [0, \alpha)\to[0,\infty)$ then $\bdry B_\varepsilon(r)$ has filling radius bounded below by a function $r'=r'(\varrho, r)$. If $\varrho$ is linear as in our case then $r' \geq \delta r$. Their result can thus be used to obtain \thmref{Theorem:lower-bounds-fillrad} in the special case that $\spt T$ is isometric to a Riemannian manifold.
Our proof for the general case uses slicing for integral currents instead of smoothing and the Lipschitz extension theorems from \secref{Section:Lip-extension} instead of the extension argument used in \cite{Greene-Petersen}. 

\begin{proof}
 Throughout the proof we view $Z$ as a subset of $X:= l^\infty(Z)$.
 Let $$0<r<\min\{2^{-(m+6)}\lambda^{-(m+1)}r_0, d(z, \spt\bdry T)\}$$ be such that $T':=\bdry(T\rstr \bar{B}(z, r))$ is an integral current with support in $S(z,r)\cap \spt T$, where $S(z, r)$ denotes the metric sphere. By the slicing theorem \cite[Theorems 5.6 and 5.7]{Ambr-Kirch-curr} this is the case for almost every $r$. Suppose, by contradiction, that $\fillrad_\infty(T')<Ar$ where $A:= \frac{1}{8(2\lambda)^{m+1}}$. We will first construct an integral current $S'$ with $\bdry S'= T'$ as follows. If $T'=0$ then set $S':= 0$. 
 If $T'\not = 0$ then choose $S\in\intcurr_m(X)$ with $\bdry S = T'$ and $\spt S\subset U_{\varepsilon}(S(z, r))\subset X$ for some
 $\varepsilon<Ar$.
Note that $\spt T'$ is a nonempty closed subset of $\spt T$ of Nagata dimension at most $m$. Furthermore, $\spt T$ is Lipschitz $m$-connected in the small and $\bar{B}(z, r_0)\cap \spt T$ is $\lambda$-linearly locally $m$-connected in $\spt T$. Therefore, by \propref{proposition:weakly-Lip-connected}, $\bar{B}(z, r_0)\cap \spt T$ is $2\lambda$-linearly locally weakly Lipschitz $m$-connected in $\spt T$. Let $f: \spt T' \to \spt T$ be the inclusion map. By \thmref{Theorem:Lip-extension-important}, there thus exists a Lipschitz extension $\bar{f}: U_{\varepsilon}(S(z, r))\to \spt T$ of $f$ with image in $U_{8(2\lambda)^{m+1}\varepsilon}(S(z, r))\cap \spt T$. Finally, set $S':= \bar{f}_{\#} S$ and note that $\bdry S' = T'$. This completes the construction of $S'$. It is now clear that in either case
 \begin{equation}\label{equation:critical-radius}
   \spt S'\subset U_{8(2\lambda)^{m+1}\varepsilon}(S(z, r))\cap \spt T.
 \end{equation}
Since $8(2\lambda)^{m+1}\varepsilon<r$ and since $z\in\spt T$, it thus follows with \eqref{equation:critical-radius} that $T'':= T\rstr \bar{B}(z,r) - S' \not= 0$; finally we have $\bdry T''=0$ and $\spt T''\subset \bar{B}(z, 2r)\cap \spt T$. 
Now, set $C:= \bar{B}(z, 2r)\cap\spt T$, endow $X':= [0, 2r]\times C$ with the product metric and define a subset of $X'$ by $Z':= \{0, 2r\}\times C$. Let $f': Z'\to \spt T$ be the $1$-Lipschitz map given by $f'(z', 0):= z'$ and $f'(z', 2r):= z$ for all $z'$. Since $2r<Ar_0$ there exists, as above, a Lipschitz extension $\bar{f}': X'\to \spt T$ of $f$. Finally, define $g(t, z'):= \bar{f}'(2rt, z')$ for all $t\in[0,1]$ and $z'\in C$ and 
set $W:= g_{\#}([0,1]\times T')$, where $[0,1]\times T'$ is the product of currents defined in Section 2.3 of \cite{Wenger-GAFA}. Note that $W\in\intcurr_{m+1}(\spt T)$ and that $\bdry W = -T''$. However, since $\hm^{m+1}(\spt T)=0$ it follows that $W=0$ and thus $T''=0$, which gives a contradiction. This completes the proof.
\end{proof}

Now, for a compact metric space $Y$ and $\varepsilon>0$, denote by $N(Y, \varepsilon)$ the smallest number of $\varepsilon$-balls needed to cover $Y$.

\bc\label{Corollary:covering-function}
 Let $Z$ be a complete metric space and $m\geq 1$. Let $T\in\intcurr_m(Z)$ with $\bdry T=0$ and suppose $\spt T$ is compact, Lipschitz $m$-connected in the small, and has Hausdorff and Nagata dimension $< m+1$. If $\spt T$ is $\lambda$-linearly locally $m$-connected up to scale $r_0$ for some $\lambda\geq 1$ and $r_0>0$ then
  \begin{equation*}
  N(\spt T, \varepsilon)\leq \frac{\lambda^{m(m+1)}\mass{T}}{C\varepsilon^m}\quad\text{for all $\varepsilon\in\left(0,2^{-(m+6)}\lambda^{-(m+1)}r_0\right)$},
 \end{equation*}
where $C$ only depends on $m$.
\ec

\begin{proof}
 Fix $0<\varepsilon < 2^{-(m+6)}\lambda^{-(m+1)}r_0$ and let $\{x_1, \dots, x_k\}\subset \spt T$ be an $\varepsilon$-separated net, that is, $d(x_i, x_j)>\varepsilon$ for all $i\not= j$. Then the balls $\bar{B}(x_j, \varepsilon/2)$ are pairwise disjoint and, by \thmref{Theorem:lower-bounds-fillrad}, we have 
\begin{equation*}
  \|T\|(\bar{B}(x_j,\varepsilon/2))\geq\frac{C'\varepsilon^m}{[16(2\lambda)^{m+1}]^m}
\end{equation*}
for every $j$ and consequently
\begin{equation}\label{equation:upper-bound-number-balls}
 k \leq \frac{[16(2\lambda)^{m+1}]^m\mass{T}}{C'\varepsilon^m}.
\end{equation}
The statement follows.
\end{proof}

\section{Proof of the main results}\label{Section:main-results}

For the proof of our main result we will need the following lemma. 

\bl\label{lemma:seq-slices}
Let $Z$ be a complete metric space and $(T_n)\subset\intcurr_m(Z)$ a bounded sequence weakly converging to some $T\in\intcurr_m(Z)$. Let $z\in Z$ and let $(z_n)\subset Z$ satisfy $z_n\to z$. Then for almost every $r>0$ there exists a subsequence $T_{n_j}$ such that $T_{n_j}\rstr \bar{B}(z_{n_j}, r)$ is a bounded sequence of integral currents converging weakly to $T\rstr \bar{B}(z, r)$, which is an integral current.
\el

The proof is a simple variation of the proof of Lemma 8.4 in \cite{Ambr-Kirch-curr}, see also Proposition 6.6 in \cite{Lang-loc-currents}.

\begin{proof}[Proof of \thmref{theorem:main}]
 Let $T$, $T_n$, and $z$, $z_n$ be as in the statement of the theorem. In the following, we view $Z$ as a subset of $l^\infty(Z)$. By \lemref{lemma:seq-slices} and \thmref{Theorem:lower-bounds-fillrad}, for almost every $0<s< 2^{-(m+6)}\lambda^{-(m+1)}r$ there exists a subsequence $T_{n_j}$ such that $T_{n_j}\rstr \bar{B}(z_{n_j}, s)$ is a bounded sequence of integral currents converging weakly to the integral current $T\rstr \bar{B}(z, s)$ and
\begin{equation*}
 \fillvol_\infty(\bdry (T_{n_j}\rstr \bar{B}(z_{n_j}, s))) \geq \frac{C's^m}{[8(2\lambda)^{m+1}]^m}
\end{equation*}
for all $j$. Since $\bdry(T_{n_j}\rstr \bar{B}(z_{n_j}, s))$ converges weakly to $\bdry(T\rstr \bar{B}(z, s))$ and since $l^\infty(Z)$ admits local cone type inequalities, Theorem 1.4 in \cite{Wenger-flatconv} implies that
\begin{equation*}
\fillvol_\infty(\bdry (T_{n_j}\rstr \bar{B}(z_{n_j}, s)) - \bdry (T\rstr \bar{B}(z, s)))\to 0
\end{equation*}
and thus 
\begin{equation*}
\|T\|(\bar{B}(z,s))\geq \fillvol_\infty(\bdry(T\rstr \bar{B}(z, s)))\geq \frac{C's^m}{[8(2\lambda)^{m+1}]^m}.
\end{equation*}
Now, \eqref{equation:vol-growth-main-thm} easily follows and hence
that $$\spt T = \{z\in Z: \text{$\|T\|$ has strictly positive $m$-dimensional lower density at $z$}\}.$$
The countable $\hm^m$-rectifiability of the set on the right-hand side was proved in Theorem 4.6 in \cite{Ambr-Kirch-curr}.
\end{proof}

\br{\rm
 The above proof in particular shows the following:  Let $(T_n)\subset\intcurr_m(Z)$ be a bounded sequence weakly converging to some $T\in\intcurr_m(Z)$ and suppose $z\in Z$ and $z_n\in\spt T_n$ are such that $z_n\to z$ for all $n$. If there exists $\delta>0$ such that each $T_n$ satisfies 
 \begin{equation}\label{equation:fillvol-good-outcome}
  \fillvol_\infty(\bdry (T_n\rstr \bar{B}(z_{n}, s))) \geq \delta s^m
 \end{equation}
  for almost all $s\in(0,\varepsilon)$ then $z\in\spt T$ and $\|T\|(\bar{B}(z,s))\geq \delta s^m$ for all $s\in(0,\varepsilon)$. It would be desirable to find weaker conditions than the $\lambda$-linear local $m$-connectedness conditions we impose in this article which still guarantee that \eqref{equation:fillvol-good-outcome} holds.
}
\er

The following example shows that if in \thmref{theorem:main} one only requires that each $\spt T_n$ has Nagata and Hausdorff dimension $\leq m+1$ instead of $< m+1$, then the theorem becomes false.

\begin{example}\label{Example:Nagata-counter}{\rm
Let $(T_n)\subset\intcurr_m(\R^{m+1})$ be a sequence with $\bdry T_n=0$ and $\spt T_n = \bar{B}(0,1)$, the closed unit ball in $\R^{m+1}$, for every $n$ and such that $\mass{T_n}\to 0$. Such a sequence can easily be obtained as a suitable infinite sum of $m$-spheres in $B(0,1)$.)
It follows that $T_n$ satisfies all the assumptions in the theorem except that $\spt T_n$ has Hausdorff and Nagata dimension $m+1$.
Furthermore, for any $z\in B(0,1)$ the constant sequence $z_n=z$ and $\lambda=1$ and $r=1-|z|$ satisfy the properties of the theorem. However, the weak limit $T$ of $T_n$ is clearly $0$.}
\end{example}

It is not clear whether \thmref{theorem:main} remains true if we only require that $\spt T_n$ has Nagata dimension at most $m$. 

Next, we turn to \thmref{Theorem:spt-rect}, the proof of which uses the following lemma.

\bl\label{lemma:upper-bound-fillvol}
 Let $Z$ be a complete metric space, $m\geq 1$, and $\varepsilon>0$. Suppose $T\in\intcurr_m(Z)$ satisfies $\bdry T= 0$ and $\fillvol_Z(T)<\varepsilon$. Then for every $z\in\spt T$ and all $r\geq 0$, $\delta>0$ there exists a set $\Omega\subset (r,r+\delta)$ of strictly positive Lebesgue measure such that $\bdry (T\rstr \bar{B}(z, s))\in\intcurr_{m-1}(Z)$ and 
 \begin{equation*}
  \fillvol_Z(\bdry (T\rstr \bar{B}(z, s)))\leq\frac{\varepsilon}{\delta}
 \end{equation*}
 for all $s\in\Omega$.
\el

\begin{proof}
Choose $S\in\intcurr_m(Z)$ with $\bdry S=T$ and $\mass{S}<\varepsilon$. Define $\varrho(y):= d(z, y)$ and note that, by the slicing theorem \cite[Theorem 5.6]{Ambr-Kirch-curr}, 
 \begin{equation*}
  \int_0^\infty\mass{\langle S, \varrho, t\rangle}\;dt\leq \mass{S}.
 \end{equation*}
 There thus exists a subset $\Omega\subset(r, r+\delta)$ of strictly positive measure, such that $\langle S, \varrho, t\rangle\in\intcurr_m(Z)$ and $\mass{\langle S, \varrho, t\rangle}\leq \varepsilon/\delta$ for all $t\in\Omega$,
   since otherwise
 \begin{equation*}
  \int_0^\infty\mass{\langle S, \varrho, t\rangle}\,dt\geq  \int_r^{r+\delta}\mass{\langle S, \varrho, t\rangle}\,dt > \varepsilon,
 \end{equation*}
 a contradiction. Now, since $\bdry(\langle S, \varrho, t\rangle) = -\langle \bdry S, \varrho, t\rangle = - \bdry (T\rstr \bar{B}(z, t))$, the proof is complete.
\end{proof}

We are ready for the proof of \thmref{Theorem:spt-rect}.

\begin{proof}[Proof of \thmref{Theorem:spt-rect}]
Let $z\in\spt T$. By \lemref{lemma:incl-spt-limit}, there exists a sequence $(z_n)$ converging to $z$ and such that $z_n\in\spt T_n$ for all $n$. It follows from \thmref{theorem:main} that \eqref{equation:growth-spt} holds and in particular that $\spt T$ is compact. 

Next, view $Z$ as an isometric subset of $l^\infty(Z)$.
We show that for every $\varepsilon>0$ the $\spt T_n$ eventually lie in the $\varepsilon$-neighborhood of $\spt T$. For this, set $A:=\spt T$ and let $A(\varepsilon):= \{z\in Z: d(z, A)\leq \varepsilon\}$. By Proposition 8.3 in \cite{Ambr-Kirch-curr} there exists for almost every $\varepsilon>0$ a subsequence $T_{n_j}$ such that $T_{n_j}\rstr A(\varepsilon)$ is a bounded sequence of integral currents weakly converging to $T$. In particular, it follows that $T'_j:=T_{n_j}\rstr [A(\varepsilon)]^c$ forms a bounded sequence converging weakly to $0$. By Theorem 1.4 in \cite{Wenger-flatconv}, $\fillvol_\infty(\bdry T'_j)\to 0$ as $j\to\infty$.
In particular, there exists $S_j\in\intcurr_m(l^\infty(Z))$ with $\bdry S_j = \bdry T'_j$ and $\mass{S_j}\to 0$. By the remark following Lemma 3.4 in \cite{Wenger-GAFA}, we may assume that $\spt S_j\subset A(2\varepsilon)$ for all $j$ large enough. Note that $T''_j:= T'_j - S_j$ satisfies $\bdry T''_j=0$, has uniformly bounded mass, and weakly converges to $0$. Therefore, $\nu_j:=\fillvol_{l^\infty(Z)}(T''_j)\to 0$, again by Theorem 1.4 in \cite{Wenger-flatconv}. Suppose now that for some $j$ large enough (to be determined later) there exists $z\in\spt T_{n_j}$ with $d(z, A)\geq 3\varepsilon$. Set
\begin{equation*}
 r':= \frac{1}{2}\min\left\{\varepsilon, 2^{-(m+6)}\lambda^{-(m+1)}r\right\}.
\end{equation*}
Since $T''_j\rstr \bar{B}(z, s) = T_{n_j}\rstr \bar{B}(z, s)$ for all $s<\varepsilon$, it follows from \thmref{Theorem:lower-bounds-fillrad} that
\begin{equation*}
 \fillvol_\infty(\bdry(T''_j\rstr \bar{B}(z,s)))\geq \frac{C's^m}{[8(2\lambda)^{m+1}]^m}
\end{equation*}
for almost every $s\in(r', 2r')$. On the other hand, \lemref{lemma:upper-bound-fillvol} implies that 
$$\fillvol_\infty(\bdry(T''_j\rstr \bar{B}(z,s)))\leq \nu_j/ r'$$
 for all $s$ in a subset of $(r', 2r')$ of strictly positive measure. If $j$ is large enough this is clearly a contraction.
This shows that $\spt T_n$ eventually lies in the $3\varepsilon$-neighborhood of $\spt T$.
Finally, let $\varepsilon>0$ and let $\{z_1,\dots, z_k\}$ be a finite and $\varepsilon$-dense set in $\spt T$. By \lemref{lemma:incl-spt-limit}, there exist sequences $z_n^j\to z_j$ with $z_n^j\in \spt T_n$. This shows that $\spt T$ lies in the $3\varepsilon$-neighborhood of $\spt T_n$ for all $n$ large enough. This shows that $\spt T_n$ converges in the Hausdorff distance to $\spt T$.
\end{proof}

\br\label{Remark:abstract-sequence-limits}{\rm
We can also deduce the following result, related to \thmref{Theorem:spt-rect}. Let $(X_n)$ be a sequence of complete metric spaces, $T_n\in\intcurr_m(X_n)$ with $\bdry T_n=0$ and such that
 \begin{equation*}
  \sup_n[\mass{T_n} + \diam(\spt T_n)]<\infty.
 \end{equation*}
Suppose further that $\spt T_n$ is compact, Lipschitz $m$-connected in the small, and has Hausdorff and Nagata dimension $< m+1$  for every $n$. If for some $\lambda\geq 1$ and $r_0>0$, every $\spt T_n$ is $\lambda$-linearly locally $m$-connected up to scale $r_0$ then it follows from \corref{Corollary:covering-function} that the sequence $(\spt T_n)$ is uniformly compact.
In particular, it follows from Gromov's compactness theorem that there exists a compact metric space $Z$ and isometric embeddings $\varphi_n:\spt T_n\hookrightarrow Z$. Moreover, there exists a subsequence $T_{n_j}$ such that $\varphi_{n_j}(\spt T_{n_j})$ converges in the Hausdorff sense to some closed subset $X\subset Z$ and $\varphi_{n_j\#}T_{n_j}$ converges weakly to some $T\in\intcurr_m(Z)$. (The latter follows from the closure and compactness theorems in \cite{Ambr-Kirch-curr}.) We conclude from \lemref{lemma:incl-spt-limit} and \thmref{theorem:main} that 
$X= \spt T$, that \eqref{equation:vol-growth-main-thm} holds at every $z\in\spt T$, and thus, by \cite{Ambr-Kirch-curr}, that $\spt T$ is countably $\hm^m$-rectifiable. Note that $X$ is the Gromov-Hausdorff limit of the sequence $(\spt T_{n_j})$.}
\er

Note that in \cite{Wenger-compactness} it is shown that whenever $(X_n)$ is a sequence of complete metric spaces and $T_n\in\intcurr_m(X_n)$ is such that
 \begin{equation}\label{equation:wenger-cpt-rmk}
  \sup_n[\mass{T_n} + \mass{\bdry T_n} + \diam(\spt T_n)]<\infty
 \end{equation}
then there exists a complete metric space $Z$, isometric embeddings $\varphi_n:\spt T_n\hookrightarrow Z$, and a subsequence $T_{n_j}$ such that $\varphi_{n_j\#}T_{n_j}$ converges with respect to the flat distance to some $T\in\intcurr_m(Z)$. In this greater generality, where only \eqref{equation:wenger-cpt-rmk} is assumed, the sequence $(\spt T_n)$ of supports need not be uniformly compact of course. Moreover, we cannot hope that a subsequence of $\spt(\varphi_{n_j\#}T_{n_j})$ converges in the Hausdorff sense. If it does then in general $\spt T$ is a strict subset.

\bigskip

Next, we give the proof of \thmref{Theorem:GH-limit}.

\begin{proof}[Proof of \thmref{Theorem:GH-limit}]
For each $n$ let $T_n$ be the integral $m$-current on $X_n$ obtained by integration over $X_n$ (the same way as the current induced by an oriented Riemannian manifold). It follows that $\bdry T_n = 0$ and, from Lemma 9.2 and Theorem 9.5 in \cite{Ambr-Kirch-curr}, that $\spt T_n = X_n$ and $\mass{T_n} \leq 2^m\omega_m^{-1}\hm^m(X_n)$ for all $n$, where $\omega_m$ denotes the volume of the unit ball in $\R^m$. It thus follows from  \corref{Corollary:covering-function} that $(X_n)$ is uniformly compact.
 In particular, by Gromov's compactness theorem, there exists a compact metric space and isometric embeddings $\varphi_n: X_n\hookrightarrow Z$. By the compactness and closure theorems of Ambrosio-Kirchheim \cite{Ambr-Kirch-curr}, possibly after passing to a subsequence, $\varphi_{n\#}T_n$ weakly converges to some $T\in\intcurr_m(Z)$. 
After passing to a further subsequence we may assume that $\varphi_n(X_n)=\spt (\varphi_{n\#}T_n)$ converges in the Hausdorff distance to a closed subset $X\subset Z$. By \lemref{lemma:incl-spt-limit} we have that $\spt T\subset X$. Furthermore, if $x\in X$ then there exists $x_n\in \spt(\varphi_{\#}T_n)$ with $x_n\to x$ and thus, by \thmref{theorem:main}, we obtain that $x\in \spt T$ and $\|T\|$ has strictly positive $m$-dimensional lower density at $x$. This shows that $\spt T = X$ and that $\spt T$ is countably $\hm^m$-rectifiable and this concludes the proof. (Note that we could have used 
\thmref{Theorem:spt-rect} that $\varphi_n(X_n)=\spt (\varphi_{n\#}T_n)$ to achieve the same.)
\end{proof}

In order to deduce \corref{Corollary:Riemannian-lin-loc} from \thmref{Theorem:GH-limit} we must show that the $M_n$ have uniformly bounded diameter. For this let $x,y\in M_n$ and let $\gamma$ be a length-minimizing geodesic joining $x$ and $y$, parametrized by arc-length. Choose points $x_1, \dots, x_k\in\image(\gamma)$ with mutual distance at least
\begin{equation*}
 \varepsilon:= \frac{r}{32(2\lambda)^{m+1}}.
\end{equation*}

It follows from \corref{Corollary:covering-function} that $k$ is uniformly bounded and thus that the length of $\gamma$ is uniformly bounded.
This establishes the uniform upper bound on diameter.

\section{Examples of cancellation and collapse}\label{Section:examples-cancellation}
Let $(T_n)$ be a bounded sequence of integral currents in some metric space $Z$ weakly converging to an integral current $T$. It was shown in \lemref{lemma:incl-spt-limit} that if the sequence $(\spt T_n)$ of supports converges in the Hausdorff sense to a closed subset $X\subset Z$ then $\spt T\subset X$. 
In the following we illustrate with some simple examples that the inclusion may be strict in general.

\begin{example}\label{example:torus}{\rm
For each $n\geq 1$, let $M_n$ be the $2$-torus $S^1\times S^1_{1/n}$ in $\R^4$, where $S^1_{1/n}$ is the circle of radius $1/n$. Then for any orientation on $M_n$, the corresponding integral currents $\Lbrack M_n\Rbrack$ converge in mass to $0$, that is $\mass{\Lbrack M_n\Rbrack}\to 0$, and thus converge weakly to $0$. On the other hand,  $M_n$ converges in the Hausdorff sense to $S^1\times \{(0,0)\}$. In particular, $\spt T$ is empty whereas the Hausdorff limit of the sequence $(\spt T_n)$ is not.
}
\end{example}

Sometimes, the phenomenon appearing in the example above is called collapse (of mass). The next example shows that the limit can be $0$ even if the mass of $T_n$ is bounded away even locally.

\begin{example}\label{example:ellipsoid}{\rm
For each $n\geq 1$ let $M_n$ be the ellipsoid in $\R^3$ given by
\begin{equation*}
 M_n := \left\{(x,y,z)\in\R^3: x^2 + y^2 + nz^2=1\right\}.
\end{equation*}
For any orientation on $M_n$ the sequence of integral currents $\Lbrack M_n\Rbrack$ converges weakly to $0$. On the other hand, the sequence $(M_n)$ converges in the Hausdorff sense to the flat disc $\{(x,y,0): x^2+y^2\leq 1\}\subset\R^3$. In particular, $\spt T$ is empty whereas the Hausdorff limit of the sequence $(\spt T_n)$ is not.
}
\end{example}

We call a phenomenon such as appearing in this example cancellation. Note that in both examples above, $M_n$ does not carry the length metric but the induced metric from the ambient Euclidean space. 

A more elaborate example of a sequence with cancellation is given as follows.

\begin{example}\label{example:cylinder-gluing}{\rm
Let $M_0$ be a closed oriented Riemannian manifold and fix $n\in\N$. Choose a collection of points, 
\begin{equation*}
\{p_1,\dots, p_{N_n} \} \subset M_0
\end{equation*}
such that $d(p_i,p_k)> 3/n$ and $M_0= \bigcup_{i}B(p_i, 10/n)$.
We choose any $r_n$ such that $0< r_n \le \min\{1/n, \injrad(M_0)/2\}$, where $\injrad(M_0)$ denotes the injectivity radius of $M_0$.
We then construct a Riemannian manifold $M_n$ by gluing $M_0\setminus \bigcup_{i=1}^{N_n} B(p_i, r_n)$
to itself and smoothing the edges so that the metric on both copies of 
$M_0 \setminus \bigcup_{i=1}^{N_n} B(p_i,2r_n)$ is preserved.
 $M_n$ converges in the Gromov-Hausdorff sense to $M_0$.
Thus there exists a compact metric space $Z$ and isometric embeddings
$\varphi_n: M_n \hookrightarrow Z$ and $\varphi_0:M_0 \hookrightarrow Z$ such that
$\varphi_n(M_n)$ Hasudorff converges to $\varphi_0(M_0)$.
However, the sequence $\varphi_{n\#}\Lbrack M_n\Rbrack$ of integral $m$-currents in $Z$ is easily seen to converge weakly to $0$. Note that there are no sequences of points with uniform
bounds on their $\lambda$ local contractibility radius due to the increasingly
dense topology on the $M_n$.
}
\end{example}

It is of course not difficult to construct sequences that exhibit partial collapse or partial cancellation.

\section{Weak convergence and Ricci curvature}\label{Section:Ricci}

Let $(M_n)$ be a sequence of closed oriented Riemannian manifolds of dimension $m$ with $\Ricci_{M_j} \geq 0$ and $\diam(M_j) \leq D_0$
 for all $n$. Denote by $T_n$ the integral $m$-current in $M_n$ induced by integration over $M_n$. Using the Bishop-Gromov volume comparison theorem it is not difficult to show that the sequence $(M_n)$ is uniformly compact and thus, by Gromov's compactness theorem, a subsequence $(M_{n_j})$ converges in the Gromov-Hausdorff sense to a compact metric space $X$. More precisely, there exists a compact metric space $Z$ and isometric embeddings $\varphi_n: M_n\hookrightarrow Z$; as in \remref{Remark:abstract-sequence-limits} it follows that for some  subsequence, which we denote again by $(M_n)$,  $\varphi_n(M_n)$ converges in the Hausdorff distance to a compact subset $X\subset Z$ and $\varphi_{n\#}T_n$ converges weakly to some $T\in\intcurr_m(Z)$.

\bt \label{Theorem:conv-Ricci}
Let $M_n$, $T_n$, $Z$, and $\varphi_n$ be as above. 
If $\Vol(M_n) \geq V_0>0$ for all $n$ then $\spt T$ coincides with the (Gromov-)Hausdorff limit $X$ of the sequence $(\varphi_n(M_n))$. Moreover,
$\|T\|$ has strictly positive $m$-dimensional lower density at every point $z\in\spt T$ and, in particular, $X$ is countably $\hm^m$-rectifiable.
\et

Under the assumptions in the theorem, there do not exist, in general,  $\lambda\geq 1$ and $r>0$ such that $M_n$ is $\lambda$-linearly locally $m$-connected up to scale $r$ for all $n$, see \cite{Perelman-example}. As the example in \cite{Menguy-inf-top-type} shows, it is not even true that for every $z\in X$ there exist $\lambda\geq 1$, $r>0$, and a sequence $z_n\in M_n$ with $\varphi_n(z_n)\to z$ such that $B(z_n, r)$ is $\lambda$-linearly locally $m$-connected in $M_n$ for $n$ large enough. In particular, \thmref{Theorem:conv-Ricci} does not come as a direct consequence of either \thmref{Theorem:spt-rect} or \thmref{theorem:main}. It should also be noted that Cheeger-Colding \cite{ChCo-PartIII} prove much stronger metric structure properties for the Gromov-Hausdorff limit $X$ than the countable $\hm^m$-rectifiability we exhibit here. Our proof uses results from Cheeger-Colding's \cite{ChCo-PartI} but not from \cite{ChCo-PartIII}.

\begin{proof} 
 Denote by $\regset\subset X$ the subset of regular points in $X$ in the sense of Cheeger-Colding, see \cite[Definition 0.1]{ChCo-PartI}.  Since $\Vol(M_n) \geq V_0>0$ for all $n$, the sequence $(M_n)$ is non-collapsed by the Bishop-Gromov volume comparison theorem. It thus follows from  Theorem 5.11 in \cite{ChCo-PartI} that at each $z\in \regset$, every tangent cone at $z$ is isometric to $\R^m$, that is, $z$ is an $m$-regular point in the sense of Cheeger-Colding.  We claim that for every $z\in\regset$ there exist $\lambda\geq 1$, $r_0>0$, and a sequence $z_n\in M_n$ with $\varphi_n(z_n)\to z$ and such that $B(z_n, r_0)$ is $\lambda$-linearly locally $m$-connected in $M_n$ for $n$ large enough. In order to do so, fix $z\in\regset$ and let $\varepsilon>0$ be sufficiently small (to be determined later). Given $\delta>0$ there exists $r>0$ such that
 \begin{equation*}
d_{GH}(B_{\R^m}(0,r), B_X(z,r)) < \delta r/2,
\end{equation*}
where $B_{\R^m}(0,r)$ and $B_X(z,r)$ denote balls in $\R^m$ and $X$, respectively.
Moreover, there is a sequence $(z_n)$ with $z_n\in M_n$ and such that $\varphi_n(z_n)\to z$ and
\begin{equation*}
d_{GH}(B_{M_n}(z_n,r), B_X(z,r)) < \delta r/2
\end{equation*}
for $n$ large enough. It follows that
\begin{equation*}
d_{GH}(B_{\R^m}(0,r), B_{M_n}(z_n,r)) < \delta r
\end{equation*}
for $n$ large enough. It now follows from Colding's volume convergence theorem \cite{Colding-volume}, more precisely from Corollary 2.19 in \cite{Colding-volume}, that if $\delta$ was chosen sufficiently small then
\begin{equation*}
 \Vol(B_{M_n}(z_n), r) \geq (1-\varepsilon) \omega_m r^m
\end{equation*} 
for $n$ large enough, where $\omega_m$ is the volume of the unit ball in $\R^m$.
Finally, if $\varepsilon>0$ was chosen small enough, Perelman's local contractibility theorem \cite{Perelman-max-vol}, see the Main Lemma therein and the remark following it, shows that $B_{M_n}(z,r_0)$ is $2$-linearly locally contractible in $M_n$ for for all $n$ large enough, where $r_0:= \nu r$ for some absolute constant $\nu\in(0,1)$. This proves our claim. It thus follows from our main result, \thmref{theorem:main}, that $\regset\subset S_T$, where $S_T$ is the set of points at which $\|T\|$ has strictly positive $m$-dimensional lower density. On the other hand, by \cite[Theorems 5.9 and 2.1]{ChCo-PartI}, we have $\hm^m(X\backslash\regset)=0$ and, in particular, $\regset$ is dense in $X$. We conclude that $\spt T = X$. 
 
 It remains to show that $\|T\|$ has strictly positive lower $m$-dimensional density at each point $z\in X\backslash\regset$. For this, we note first that, by Lemma 9.2 and Theorem 9.5 in \cite{Ambr-Kirch-curr}  and the fact that $\hm^m(X\backslash\regset)=0$,  we have
\begin{equation}\label{equation:mass-growth}
 \|T\|(B(z,r))\geq m^{-m/2}\hm^m(S_T\cap B(z, r)) = m^{-m/2}\hm^m(X\cap B(z, r))
\end{equation}
for all $r>0$. On the other hand, by Theorem 5.9 in \cite{ChCo-PartI} and the Bishop-Gromov volume comparison theorem,
\begin{equation}\label{equation:hm-growth}
 \hm^m(X\cap B(z,r)) \ge \frac{V_0r^m}{(\diam X)^m}
\end{equation}
for all $0<r<\diam X$.
From \eqref{equation:mass-growth} and \eqref{equation:hm-growth} it thus follows that $\|T\|$ has strictly positive $m$-dimensional lower density at $z$. This concludes the proof of the last statement.
\end{proof}


We end this section by noting that if the Ricci curvature condition is replaced by a scalar curvature
condition then \thmref{Theorem:conv-Ricci} no longer holds even when the sequence of manifolds
is known to converge in the Gromov-Hausdorff sense.  There can be cancellation without
collapse, as the following example shows.

\begin{example} \label{example-scalar}{\rm
Let $M_0$ be a closed oriented $m$-dimensional Riemannian manifold with positive scalar curvature and $m\ge 3$. Construct $M_n$
 analogously to the $M_n$ in \exref{example:cylinder-gluing}, however using the Gromov-Lawson gluing construction \cite{Gromov-Lawson} on tiny balls very close to the $p_i$
(where $M_0$ is sufficiently flat). It can be achieved that $M_n$ has positive scalar curvature and the metric on $M_0\setminus \bigcup_{i=1}^{N_n} B(p_i, 2r_n)$ is preserved.  
Then these $M_n$ still converge in the Gromov-Hausdorff sense to $M_0$ and  the
$\varphi_{n\#}\Lbrack M_n\Rbrack$ converge weakly to $0$.
}
\end{example}


\bigskip

\bigskip

\appendix
\section{Limits of uniformly geometric contractible sequences of manifolds\\*[4ex] { BY RAANAN SCHUL AND STEFAN WENGER}\\*[5ex]}

The purpose of this appendix is to prove the following theorem which shows that $\lambda$-linear local $m$-connectedness cannot be replaced by geometric contractibility in \corref{Corollary:Riemannian-lin-loc}.

\bt\label{theorem:non-rect}
 For every $m\geq 2$ and every $\alpha\in(0,1)$ there exist $C, r>0$ and a sequence of Riemannian metrics $g_n$ on the $m$-sphere $S^m$ such that $M_n:= (S^m, g_n)$ has
 \begin{equation*}
  \sup_{n}\left[\Vol(M_n) + \diam M_n\right]<\infty,
 \end{equation*}
 and satisfies the following properties:
 Each $M_n$ has local geometric contractibility function $\varrho(s):= Cs^\alpha$, $s\in[0,r)$, and $M_n$ converges in the Gromov-Hausdorff sense to a compact metric space $X$, homeomorphic to $S^m$, which satisfies $0<\hm^m(X)<\infty$ but which is not countably $\hm^m$-rectifiable.
\et

In the proof we will need the following:
Let $N, m\geq 2$ and set $a:= N^{-1/m}$. Define a metric $d_\infty$ on $Z:= \{1,\dots, N\}^{\N}$ such that for given $(z_l), (z'_l)\in Z$ we have $d_\infty((z_l), (z'_l)):= a^j,$ where $j$ is the smallest index, or $0$, for which $z_j\not=z'_j$. Note that $d_\infty$ is an ultrametric, that is,
\begin{equation}\label{equation:ultrametric}
 d_\infty(z, z') \leq \max\left\{d_\infty(z,z''), d_\infty(z',z'')\right\}
\end{equation}
for all $z,z',z''\in Z$.

\bl\label{lemma-appendix}
 The metric space $(Z, d_\infty)$ satisfies $0<\hm^m(Z)<\infty$, but is not countably $\hm^m$-rectifiable.
\el

\begin{proof}
 We first show that $Z$ has finite $\hm^m$-measure. For this, let $\delta>0$ and choose $n_0$ large enough so that $2a^{n_0}<\delta$. Since $Z$ can be covered by $N^{n_0}$ closed balls of radius $a^{n_0+1}$ we obtain
\begin{equation*}
 \hm^m_\delta(Z) \leq N^{n_0} \omega_ma^{(n_0+1)m}= \omega_ma^m, 
\end{equation*}
where $\omega_m$ is the volume of the unit ball in $\R^m$.
This shows that $Z$ has finite $\hm^m$-measure.
Next, we construct a Frostman measure on $Z$ and use it to show that $\hm^m(Z)>0$. Let $\pi_n: Z\to \{1,\dots, N\}$ be the projection onto the $n$-th coordinate. Let $\mu$ be the unique Borel probability measure on $Z$ such that for all $n\geq 1$ and all $z_1, \dots, z_n\in\{1,\dots, N\}$
\begin{equation*}
 \mu(\{\pi_i = z_i: i=1, \dots, n\}) = \frac{1}{N^n}.
\end{equation*}
We claim that for all $z\in Z$ and all $r\in (0,a)$ we have
$\mu(\bar{B}(z,r))\leq a^{-2m}r^m$. Indeed, let $n\in\N$ be such that $a^{n+1}<r\leq a^n$. Then every $z'\in \bar{B}(z, r)$ satisfies $\pi_i(z') = \pi_i(z)$ for $i=1, \dots, n-1$ 
and hence $\mu(\bar{B}(z,r))\leq N^{-n+1} \leq a^{-2m}r^m$. Now, let $\varepsilon>0$ and let $(B_i)$ be a covering of $Z$ by sets of diameter less than $a/2$ and such that
\begin{equation*}
 \sum_{i=1}^\infty \omega_m\left(\frac{\diam B_i}{2}\right)^m\leq \hm^m(Z) + \varepsilon.
\end{equation*}
For each $i$ set $r_i:= \diam B_i$ and choose $z_i\in B_i$ arbitrary. Then
\begin{equation*}
 \sum_{i=1}^\infty r_i^m \geq  2^{-m}a^{2m}\sum_{i=1}^\infty \mu(\bar{B}(z_i, r_i))\geq 2^{-m}a^{2m}\mu(Z) = 2^{-m}a^{2m}.
\end{equation*}
Since $\varepsilon>0$ was arbitrary, it follows that $\hm^m(Z)\geq 2^{-m}a^{2m}\omega_m>0$.

We finally show that for any $K\subset\R^m$ Borel and any Lipschitz map $\varphi : K\to Z$ we have $\hm^m(\varphi(K))=0$. Indeed, if this is not the case, then it follows from \cite{Kirchheim} that there exists $K'\subset K$ of strictly positive Lebesgue measure, a norm $\|\cdot\|$ on $\R^m$, and $x\in K'$ a Lebesgue density point of $K'$ such that the map $\varphi\on{K'}$ is an approximate isometry around $x$ when viewed as a map from $(K', \|\cdot\|)$ to $Z$. Now, for any $\varepsilon>0$ there exist $x', x''\in K'$ with $\|x-x'\|<\varepsilon$ and $\|x-x''\|\leq \frac{3}{4}\|x-x'\|$ and $\|x'-x''\|\leq \frac{3}{4}\|x-x'\|$. On the other hand, if $k$ is such that  $d_\infty(\varphi(x), \varphi(x'))=a^k$ then, by \eqref{equation:ultrametric}, either $d_\infty(\varphi(x), \varphi(x''))\geq a^{k}$ or $d_\infty(\varphi(x'), \varphi(x''))\geq a^{k}$. Since $\varphi\on{K'}$ is an approximate isometry around $x$, this yields a contradiction. This completes the proof.
\end{proof}

We use the above lemma to prove the theorem.

\begin{proof}[Proof of \thmref{theorem:non-rect}]
 Let $m, L\geq 2$ and $\alpha\in(0,1)$, set $a:= L^{-1}$ and $N:= L^m$, and let $\lambda\in(0,a)$ such that $\alpha\log \lambda = \log a$. For each $n\geq 1$, let $P_n$ be the boundary of the $(m+1)$-dimensional rectangle  $[0,\lambda^n]^m\times [0,a^n]$. We call $[0,\lambda^n]^m\times \{0\}$ the base of $P_n$ and $[0,\lambda^n]^m\times \{a^n\}$ the roof of $P_n$.
Let $X_0:=[0,L\lambda]^m$ have the Euclidean metric, divide $X_0$ into $N$ equal cubes of edge length $\lambda$ each, and glue to each cube a copy of $P_1$ along its base (and `delete' the interior of the base after gluing). Endow the so obtained space with the length metric and call it $X_1$. We call the copies of $P_1$ in $X_1$ the $1$-towers of $X_1$. Suppose now that we have constructed $X_n$ for some $n\geq 1$. In order to construct $X_{n+1}$ from $X_n$, center a cube of edge length $L\lambda^{n+1}$ on the roof of each $n$-tower in $X_n$, note for this that $L\lambda^{n+1} < \lambda^n$, then divide the cube into $N$ equal cubes of edge length $\lambda^{n+1}$ and glue a copy of $P_{n+1}$ along its base to each cube (and `delete' the interior of the base after gluing). Endow the new space with the length metric and call it $X_{n+1}$. The copies of $P_{n+1}$ in $X_{n+1}$ are called $(n+1)$-towers in $X_{n+1}$. Figure \ref{figure} below shows how towers stand on top of each other.

\begin{figure*}[h]
 \begin{minipage}[b]{7cm}
 \psfrag{aT}{$k$-tower}
 \psfrag{bT}{{\footnotesize $(k+1)$-tower}}
  \psfrag{cT}{{\footnotesize \, $(k+1)$-tower}}
 \psfrag{dT}{{\tiny \, $(k+2)$-tower}}
 \psfrag{f}{$\lambda^k$}
 \psfrag{h}{{\footnotesize $L\lambda^{k+1}$}}
 \psfrag{r}{{\tiny $\lambda^{k+1}$}}
         \includegraphics*[width=6cm]{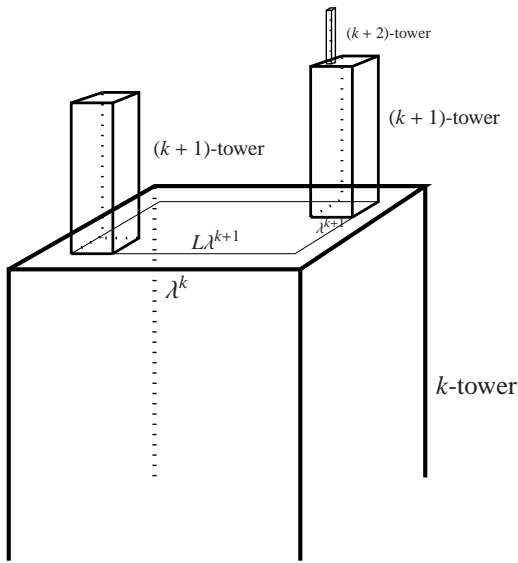}
         \caption{Some of the towers in $X_n$}\label{figure}
      \end{minipage}
\end{figure*}
Clearly, each $X_n$ is biLipschitz homeomorphic to the cube $X_0$ and satisfies
\begin{equation*}
 \hm^m(X_n)= (L\lambda)^m + 2m\sum_{k=1}^{n}\mu^k\leq (L\lambda)^m + \frac{2m}{1-\mu},
\end{equation*}
where $\mu:= aN\lambda^{m-1}<1$.
For $k<n$, the distance $d$ from any point at the base of a given $k$-tower in $X_n$ to any point on its roof satisfies $a^k \leq d \leq (m+1) a^k$.
It follows easily that all $X_n$ have uniformly bounded diameter.
Now, each $X_n$ has local geometric contractibility function  $\varrho: [0,\lambda/2)\to [0,\infty)$ given by $\varrho(s):= Cs^\alpha$, where $C = 2^\alpha(m+2-a)(1-a)^{-1}L$. Indeed, let $z\in X_n$ and $s\in(0,\lambda/2)$. Let $k\in\N$ be such that $\lambda^{k+1}\leq 2s < \lambda^k$. Suppose that $B(z,s)$ intersects a $p$-tower but no $(p-1)$-tower in $X_n$. If $p<k$ then $B(z,s)$ cannot intersect an $l$-tower with $l> p+1$ since otherwise $2s\geq a^{p+1}$ and hence $p\geq k$, a contradiction. Since $2s<\lambda^{p+1}$ it thus follows that $B(z,s)$ is contractible within itself. 
Now, if $p\geq k$ then $B(z, s)$ is contractible within the ball $B(z, s')$ where
\begin{equation*}
 s' = \lambda^k  + \frac{a^p(m+1)}{1-a}.
\end{equation*}
It is trivial to check that $s'\leq \varrho(s)$.  This shows that $\varrho$ is indeed a geometric contractibility function for $X_n$.

Let now $X_\infty$ be the Gromov-Hausdorff limit of the sequence $X_n$. We first show that $X_\infty$ has finite $\hm^m$-measure. For this, let $\delta>0$ and choose $n_0$ so large that $2(m+1)(1-a)^{-1}a^{n_0}<\delta$. Fix an $n_0$-tower in $X_\infty$ and a point $x$ at its base. Then the closed ball of radius $(m+1)(1-a)^{-1}a^{n_0}$ around $x$ contains the given $n_0$-tower and all the towers on its roofs (and the towers on their roofs, and so on). Since $X_\infty$ contains exactly $N^{n_0}$ $n_0$-towers, we see that
\begin{equation*}
 \hm_\delta^m(X_\infty) \leq \hm^m(X_{n_0-1}) + N^{n_0} \omega_m \left[\frac{(m+1)a^{n_0}}{1-a}\right]^{m}
    \leq (L\lambda)^m + \frac{2m}{1-\mu} + \omega_m \left[\frac{m+1}{1-a}\right]^{m}.
\end{equation*}
Since $\delta>0$ was arbitrary, this shows that indeed $\hm^m(X_\infty)<\infty$. Next, the metric space $Z$ from the lemma above admits a biLipschitz embedding into $X_\infty$. Since $Z$ is not countably $\hm^m$-rectifiable by \lemref{lemma-appendix}, it follows that $X_\infty$ is not countably $\hm^m$-rectifiable either.

Finally, in order to obtain a closed oriented Riemannian manifold $M_n$ we glue $X_n$ along its boundary to the boundary of an $m$-cube and then smooth it. In this way, we can also achieve that each $M_n$ is biLipschitz homeomorphic to the standard $m$-sphere $S^m$ and that the Gromov-Hausdorff limit is homeomorphic to $S^m$. This concludes the proof of our theorem.
\end{proof}

\bibliographystyle{abbrv}
\bibliography{Sormani-Schul-Wenger-bibli}

\end{document}